\documentclass[number,citesort,seceqn,dvips]{arxbj}
\usepackage{upgreek}
\usepackage[left]{vmkcol}
\usepackage{multirow}

\aid{0}
\volume{18}
\issue{1}
\pubyear{2012}
\firstpage{64}
\lastpage{99}
\doi{10.3150/10-BEJ328}

\makeatletter
\newcommand{\eqref}[1]{(\ref{#1})}
\newcolumntype{d}[1]{D{.}{.}{#1}}

\newtheorem{theorem}{Theorem}[section]
\newremark{remark}{Remark}[section]
\newtheorem{proposition}{Proposition}[section]
\newtheorem{lemma}{Lemma}[section]

\newcommand{\RR}{\mathbb{R}}
\newcommand{\NN}{\mathbb{N}}
\newcommand{\MM}{\mathbb{M}}
\newcommand{\MMM}{{\mathcal M}}
\newcommand{\AAA}{{\mathcal A}}
\newcommand{\BBB}{{\mathcal B}}
\newcommand{\DDD}{{\mathcal D}}
\newcommand{\EEE}{{\mathcal E}}
\newcommand{\FFF}{{\mathcal F}}
\newcommand{\HHH}{{\mathcal H}}
\newcommand{\LLL}{{\mathcal L}}

\makeatother

\begin{document}
\begin{frontmatter}

\title{Limit experiments of GARCH}
\runtitle{Limit experiments of GARCH}

\begin{aug}
\author[1]{\fnms{Boris} \snm{Buchmann}\thanksref{1}\ead[label=e1]{Boris.Buchmann@anu.edu.au}}
\and
\author[2]{\fnms{Gernot} \snm{M\"{u}ller}\thanksref{2}\corref{}\ead[label=e2]{mueller@ma.tum.de}}
\runauthor{B. Buchmann and G. M\"{u}ller}
\address[1]{Mathematical Sciences Institute,
School of Finance, Actuarial Studies and Applied Statistics,
Australian National University,
ACT 0200,
Australia. \printead{e1}}
\address[2]{Zentrum Mathematik,
Technische Universit\"{a}t M\"{u}nchen,
Boltzmannstra{\ss}e 3,
85748 Garching,
Germany. \printead{e2}}
\end{aug}

\received{\smonth{1} \syear{2010}}
\revised{\smonth{7} \syear{2010}}

%
\begin{abstract}
GARCH is one of the most prominent nonlinear time series models, both widely
applied and thoroughly studied. Recently, it has been shown that the
COGARCH model (which was introduced a
few years ago by Kl\"{u}ppelberg, Lindner and Maller) and Nelson's
diffusion limit are the only
functional continuous-time limits of GARCH in distribution. In contrast
to Nelson's
diffusion limit, COGARCH reproduces most of the \textit{stylized facts}
of financial time series.
Since it has been proven that Nelson's diffusion is not
asymptotically equivalent to GARCH in deficiency, in
the present paper, we investigate the relation between GARCH and
COGARCH in Le Cam's
framework of statistical equivalence.
We show that GARCH converges generically to COGARCH, even in deficiency,
provided that the volatility processes are observed. Hence, from
a theoretical point of view, COGARCH can indeed be considered as
a continuous-time equivalent to GARCH. Otherwise, when the observations
are incomplete, GARCH still has a limiting experiment, which we call MCOGARCH,
which is not equivalent, but nevertheless quite similar, to COGARCH.
In the COGARCH model, the jump times can be more random than for the
MCOGARCH, a fact
practitioners may see as an advantage of COGARCH.
\end{abstract}

%

\end{frontmatter}
%

\section{Introduction}\label{s1}
Since the seminal papers of Engle \cite{Engle82} and Bollerslev
\cite{Bo86} the discrete-time GARCH methodology has become a
widely applied tool in the modeling of heteroscedasticity in
financial times series. On the other hand, continuous-time models are
very useful, for instance, in option pricing, as shown by Black and
Scholes \cite{BS73} and Merton \cite{Me73}, in the
analysis of tick-by-tick data and for modeling irregularly spaced time
series.

In the 1990s, researchers tried to bridge the gap between continuous
and discrete time. Nelson~\cite{Ne90} showed that an
appropriately parametrized GARCH can be seen as a discrete-time
approximation of a bivariate diffusion model on an approximating time
grid. However, this diffusion model does not capture most of the
so-called \textit{stylized facts} reflecting empirical findings in
financial time series: for example, volatility exhibits heavy tails,
jumps upward and clusters on high levels. To overcome the shortcomings
of the diffusion model, Kl\"{u}ppelberg \textit{et al.} \cite{KLM04}
introduced a new continuous-time GARCH model which they called COGARCH.
In contrast to the bivariate diffusion, this model exhibits many of the
stylized facts. We refer the reader to Fasen \textit{et al.}
\cite{FKL06} for an extensive discussion of the stylized facts and the
various competing volatility models proposed in the literature.

Recently, Kallsen and Vesenmayer \cite{KV07} and Maller \textit{et
al.} \cite{MMS08} have identified COGARCH as a functional limit
of GARCH in distribution. Most notably, Kallsen and
Vesenmayer~\cite{KV07} have argued that Nelson's diffusion and COGARCH
are the only possible limits of GARCH in distribution in a
semimartingale setting.

The passage from discrete to continuous time has an obviously appealing
practical purpose: we can estimate the underlying continuous-time model
parameters by a time series formulation and plug them into the
continuous-time limit for other purposes. As argued by Wang
\cite{Wa02}, such a passage is, in general, only justified if the
corresponding statistical experiments converge in Le Cam's framework of
deficiency (see Le Cam \cite{LC86}, Le Cam and Young
\cite{LY90} and Strasser \cite{St85}).

In particular, Wang \cite{Wa02} (see also Brown \textit{et al.}
\cite{BWZ03}) showed, assuming independent Gaussian innovations, that
Nelson's diffusion approximation of GARCH is not valid in deficiency:
the innovations encounter the two models in an intrinsically different
way. Whereas GARCH is driven by one-dimensional innovations, its
diffusion limit is driven by planar Brownian motion.

In contrast to Nelson's approximation, COGARCH is driven by a L\'{e}vy
process of only one dimension, thereby mimicking one of the key
features of GARCH. Naturally, the following questions arise: Are the
approximations of COGARCH by GARCH, as proposed by Kallsen and
Vesenmayer \cite{KV07} and Maller \textit{et al.} \cite{MMS08}, also valid in
deficiency? Does the limiting model depend on the underlying sampling
scheme?

Dealing with Le Cam's distance in deficiency is a challenging task. In
particular, asymptotic equivalence results for dependent data are very
scarce; see Dalalyan and Rei{\ss} \cite{DR06} for an overview. Further
obstacles arise from the intrinsic heteroscedasticity of GARCH.
Therefore, in this paper, we restrict ourselves to compound Poisson
processes as driving L\'{e}vy process and assume that the innovations
are randomly thinned. This approximation scheme also occurs in the
papers \cite{KV07} and \cite{MMS08}. Random thinning is a standard
limiting procedure in many other areas of probability theory and
statistics. In particular, we mention the peak-over-threshold method in
extreme value theory (see Remark~\ref{KVMMS}(ii)). In contrast to our
approximation scheme, most papers on statistical equivalence deal with
aggregated innovations where the experiments are compared to Gaussian
shift experiments (see Brown and Low \cite{BL96},
Nussbaum \cite{Nu96}, Grama and Neumann \cite{GN06} and
Carter~\cite{C07} and references therein; see also Milstein and
Nussbaum \cite{MN98} with potential applications to time series
analysis). We point out once again that for the GARCH, aggregated
innovations lead to the diffusion limit investigated by
Nelson \cite{Ne90} and Wang~\cite{Wa02}.

The paper is organized as follows. Section \ref{mainresult} contains
our main results. We introduce the experiments and sampling schemes in
Section \ref{secGARCH}. In Section \ref{secIncomplete}, we construct a
limiting experiment for randomly thinned GARCH with unobserved conditional
variances. As shown in
Section~\ref{secIncompletecomparison}, using both theoretical and
numerical methods, this experiment is generically not equivalent to
COGARCH. However, if the conditional variances are observable in full,
then all experiments are generically (asymptotically) equivalent to
COGARCH -- this is shown in Section~\ref{secComplete}. We conclude in
Section \ref{summary}. In Sections \ref{proofA}--\ref{SUBSEC4.9} we
give the proofs of all theorems and propositions from
Section \ref{mainresult}. Section \ref{proofA} contains the proof of
Theorem \ref{th2}, Section \ref{proofB} the proof of Theorem \ref{th2a}
and Section \ref{proofth1} the proof of Theorem \ref{th1}. The proofs
of all propositions in Section~\ref{secComplete} are reported in
Section \ref{SUBSEC4.9}. In the \hyperref[appendix]{Appendix}, we review some of the basic
notions of Le Cam's convergence in deficiency.

\vspace*{2pt}\section{Main results}\label{mainresult}\vspace*{2pt}
\subsection{GARCH-type experiments in discrete and continuous time}\label{secGARCH}
 For all $n\in\NN$ we consider an $n$-dimensional
vector $Z_{n}=(Z_{n,k})_{1\le k\le n}$ with distribution
%
\begin{equation}\label{Znthinn}
\LLL(Z_n) =  \bigl((1 - p_n)\varepsilon_0+p_n Q_n \bigr)^{\otimes
n} ,
\end{equation}
where, for all $n\in\NN$, $p_n\in(0,1)$ and $Q_n$
is a probability measure on the Borel field $\BBB(\RR)$. Here,
$\epsilon_0$ denotes the Dirac measure with total mass in zero.

The parameter $p_n$ modulates our random thinning. In accordance with
the law of rare events, we assume that the following limit exists in
$(0,\infty)$:
%
\begin{equation}\label{lre}
\gamma = \lim_{n\to\infty}np_n\in (0,\infty) .
\end{equation}
In the sequel, we will encounter several GARCH-type processes, all
indexed by $\theta\in[0,\infty)^4$. In discrete time, processes will be
further indexed by $n\in\NN$ and a suitable parametrization. Throughout
this paper, a parametrization is a pair $(\Theta,(H_n)_{n\in\NN}),$
where $\Theta$ is a non-empty subset of $[0,\infty)^4,$ and for all
$n\in\NN$, $H_n$ is a mapping
$H_n=(h_{0,n},\beta_n,\alpha_n,\lambda_n)\dvtx\Theta\to[0,\infty)^4$. Here,
$h_0$ ($h_{0,n}(\theta)$) denotes the unknown initial value of the
volatility $h_0$, which is treated as an additional unknown parameter
in this paper. For the corresponding continuous-time limits, $\alpha$ is the mean
reversion parameter of the volatility processes and $\beta/\alpha$ the
asymptotically stable fixed point of the (unperturbated) volatility SDE
(ODE). $\lambda$~is a
scaling parameter for the corresponding jumps of the volatility
processes.

For a parametrization $(\Theta,(H_n)_{n\in\NN}),$ we consider the
sequence of partial sums corresponding to a \textit{randomly thinned}
GARCH model indexed by $\theta\in\Theta$ and $n\in\NN$, defined by
%
\begin{eqnarray}\label{G1}
G_{n}(k)&=&G_{n}(k  -  1)+h^{1/2}_n(k  -  1)
Z_{n,k} ,\qquad  G_n(0) = 0 ,\nonumber\\
h_{n}(k)&=&\beta_n(\theta)+\alpha_n(\theta)
h_{n}(k  -  1)+
\lambda_n(\theta)  h_{n}(k  -  1) Z_{n,k}^2 ,\\
h_{n}(0)&=&h_{0,n}(\theta),\qquad 1\le k\le
n ,\ \theta\in\Theta,\nonumber
\end{eqnarray}
where
$H_n(\theta)=(h_{0,n}(\theta),\beta_n(\theta),\alpha_n(\theta
),\lambda_n(\theta))$
for all $\theta\in\Theta$. Note that this specification of a~GARCH does
not quite follow the traditional one, but enumerating the indices
generates the same processes. Also, observe that the definition of
$(G_n,h_n)$ in \eqref{G1} depends on the choice of
$(\Theta,(H_n)_{n\in\NN})$.

Provided that $Q_n$ converges weakly to some probability measure $Q$,
the limit in \eqref{lre} sets up convergence in distribution of
$\sum_{k=1}^{[n\cdot]} Z_{n,k}$ to a compound Poisson process with rate
$\gamma$ and jump distribution $Q$ as $n\to\infty$. For a choice of
$(\Theta,(H_n)_{n\in\NN})$ it is thus natural to ask whether the limit
of $(G_n([nt]),h_n([nt]))_{0\le t\le1}$ in distribution exists along
$H_n(\theta)$ as $n\to\infty$ for fixed $\theta\in\Theta$. In
\cite{KV07} and \cite{MMS08}, such parametrizations have been
successfully constructed. Moreover, the corresponding continuous-time
limit equals COGARCH driven by a compound Poisson process.

COGARCH is a process $(G,h)=(G(t),h(t))_{0\le t\le1}$ that is indexed
by $\theta=(h_0,\beta,\alpha,\lambda)\in[0,\infty)^4$ and
determined as
the unique pathwise solution of the following system of integral
equations:
%
\begin{eqnarray}\label{COG}
G(t)&=&\int_{[0,t]\times\RR} h^{1/2}(s-)z N(\mathrm{d}s,\mathrm{d}z) ,\nonumber\\ [-8pt]\\ [-8pt]
h(t)&=&h_0+ \int_{[0,t]} \beta - \alpha h(s-)\,\mathrm{d}s+\lambda
\int_{[0,t]\times\RR}  h(s-)z^2 N(\mathrm{d}s,\mathrm{d}z),\nonumber
\end{eqnarray}
where $N$ is a Poisson point measure on $[0,1]\times\RR$ with intensity
$\gamma\ell\otimes Q$.

In the sequel, we restrict our analysis to the two following sampling
schemes:
\begin{itemize}
\item \textit{Incomplete observations}: only $G$ and $G_n$ ($n\in\NN$) are
observable in full, whereas the corresponding volatility processes, $h$
and $h_n$ ($n\in\NN$), are unobservable.
\item \textit{Complete
observations}: both processes $(G,h)$ and $(G_n,h_n)$ are observable in
full.
\end{itemize}
We will deal with these two sampling schemes in
Sections \ref{secIncomplete}--\ref{secIncompletecomparison} and
Section \ref{secComplete}, respectively. Not surprisingly, a simpler
theory is in place in the case of complete observations. In a more
realistic scenario, where observations of the volatility processes are
not available, results are more difficult due to the nonlinearity of
(CO)GARCH.

Throughout the paper, $ $ the space of right-continuous functions
$g\dvtx[0,1]\to\RR^d$ with left limits on $[0,1]$ is denoted by $D_d$. We
endow $D_d$ with the $\sigma$-algebra $\DDD_d$, generated by point
evaluations (see Billingsley \cite{Bi68}). Furthermore, let
$\MM_d$ be the space of all non-negative point measures on
$[0,1]\times\RR^d$ with finite support. We equip this space with the
$\sigma$-algebra $\MMM_d$ generated by the point evaluations
$A\mapsto
\mu(A)$, $A\in\BBB([0,1]\times\RR^d)$, $\mu\in\MM_0$ (see
Reiss \cite{Re93}, pages 5--6).

The trace of the Borel field in
$\overline{\RR}^d=(\RR\cup\{-\infty,\infty\})^d$ with respect to
$A\subseteq\overline{\RR}^d$ is denoted by $\BBB(A)$. The Lebesgue
measure on $\BBB(\RR)$ and the Dirac measure with total mass in $x$ are
denoted by $\ell$ and $\epsilon_x$, respectively. If $(E,\AAA)$ is a
measurable space and $X$ is a random element taking values in
$(E,\AAA),$ then its distribution is denoted by $\LLL(X)$. Whenever
this distribution depends on a parameter $\theta,$ we employ the
notation $\LLL_\theta(X)$. If $(E_i,\AAA_i)$, $i=1,2$, are measurable
spaces and $X\dvtx E_1\to E_2$ is $\AAA_1/\AAA_2$ measurable, then $\mu^X$
denotes the image of a measure $\mu$ under $X$.

We refer to the \hyperref[appendix]{Appendix} and \cite{St85} for unexplained notation
relating to convergence in deficiency.

\subsection{Limit experiments of GARCH (incomplete
observations)}\label{secIncomplete}\vspace*{-3pt}
 In this subsection, we assume that the
vola\-tility processes are unobservable. To pursue our program, we
introduce another class of processes. Let $(\widehat G, \hat
h)=(\widehat G(t), \hat h(t))_{0\le t\le1}$ be the unique pathwise
solution of the following system of integral equations:
%
\begin{eqnarray}\label{hatEEE}
\widehat G(t)&=&\int_{[0,t]\times\RR}
\hat h^{1/2}(s-)z N(\mathrm{d}s,\mathrm{d}z) ,\nonumber\\ [-9pt]\\ [-9pt]
\hat h(t)&=&h_0+\int_{[0,t]} \beta - \alpha\hat
h(s-)\,\mathrm{d}T_{N} (s)+\lambda\int_{[0,t]\times\RR}  \hat
h(s-)z^2
N(\mathrm{d}s,\mathrm{d}z) ,\nonumber
\end{eqnarray}
where $\theta=(h_0,\beta,\alpha,\lambda)\in[0,\infty)^4$. Here,
$T\dvtx\MM_1\to D_1$, $\sigma\mapsto T_\sigma$ is defined as follows: If,
for some $m\in\NN$, $0=t_0<t_1<\cdots<t_{m}< 1$ and $x_1,\dots,x_m\in
\RR$,
$\sigma\in\MM$ admits a representation of the form
$\sigma=\sum_{k=1}^m\varepsilon_{(t_k,x_k)}$, where
$0=t_0<t_1<\cdots<t_{m}< 1$, then we set
\begin{eqnarray}\label{timechange}
T_{\sigma}(t)&=&
\frac{t - t_k}{m(t_k - t_{k - 1})}+\frac km ,\qquad  t\in
[t_{k-1},t_k),  1\le k\le m ,\nonumber\\ [-9pt]\\ [-9pt]
T_{\sigma}(t)&=&\frac{t - t_m}{m(t_m - t_{m - 1})}+1,\qquad
t\in[t_{m},1].\nonumber
\end{eqnarray}
If such a representation does not exist, then we set $T_{\sigma}(t)=t$
for all $t\in[0,1]$.

We call $(\widehat G,\hat h)$ the MCOGARCH, an acronym referring to
\textit{modified} COGARCH. To illustrate the difference between COGARCH
and MCOGARCH, we next consider a simpler representation of $\widehat G$
(we will return to \eqref{hatEEE} in our analysis in
Section \ref{secComplete}).

To this end, let $\nu=(\nu(t))_{0\le t\le1}$ be a Poisson process with
rate $\gamma$ and $(Z_k)_{k\in\NN}$ be a~sequence of independent random
variables, independent of $\nu$. By solving the integral equations for
$\hat h$ in~\eqref{hatEEE}, we observe that
%
\begin{equation}\label{defhatG}
\LLL_\theta(\widehat
G) = \LLL_{\theta} \Biggl(\sum_{k=1}^{\nu(\cdot)} \hat
h^{1/2}_{\nu(1),k,\theta} Z_k \Biggr)  ,
\end{equation}
where, for
$k,m\in\NN$, $k\ge2$, we set
%
\begin{eqnarray}\label{defhhat}
\hat
h_{m,k,\theta}&=&\frac{\beta}{\alpha}  (1 - \mathrm{e}^{-\alpha
/m} )+\mathrm{e}^{-\alpha/m}\hat h_{m,k-1,\theta} [1+\lambda Z_{k-1}^2] ,\nonumber\\ [-9pt]\\ [-9pt]
\hat h_{m,1,\theta}&=&\frac{\beta}{\alpha}  (1 - \mathrm{e}^{-\alpha
/m} )+\mathrm{e}^{-\alpha/m}h_0\nonumber
\end{eqnarray}
for $\theta=(h_0,\beta,\alpha,\lambda)\in[0,\infty)^4$, $\alpha>0$,
with the convention that $\sum_\varnothing=0$. Here, we extend the
definition of $\hat h_{m,k,\theta}$ to
$\theta=(h_0,\beta,0,\lambda)\in[0,\infty)^4$ by taking
$\alpha\downarrow0$ in~\eqref{defhhat}.

In view of \eqref{defhhat}, note that the magnitudes of the jumps of
$\widehat G$ (in space) depend on their multiplicities and the sizes of
innovations, but not on their arrival times. This attribute is not
shared by COGARCH. Therefore, it is, to some extent, justified to
speak of $\widehat G$ and $G$ as experiments driven by two and three
sources of randomness, respectively: the number of jumps, the
innovations and the arrival times.\vadjust{\goodbreak}

As no information about the volatility processes is assumed in this
subsection, we consider the following experiment of MCOGARCH type:
%
\begin{equation}\label{MCOGARCH}
\widehat\EEE =  \bigl(D_1,\DDD_1, (\LLL_{\theta}(\widehat
G) )_{\theta\in[0,\infty)^4} \bigr) .
\end{equation}
For a parametrization $(\Theta,(H_n)_{n\in\NN}),$ we consider the
corresponding GARCH experiments in discrete time,
%
\begin{equation}\label{experiments}
\EEE_{n,H_n}(\Theta) =  (\RR^n,\BBB(\RR^n), (\LLL
_{\theta}(G_n) )_{\theta\in\Theta}),\qquad
n\in\NN ,
\end{equation}
where, for $n\in\NN$, $G_n=(G_n(k))_{1\le k\le n}$ is defined by
\eqref{G1} via the parametrization $(\Theta,(H_n)_{n\in\NN})$. We write
$\EEE_{n,H_n}=\EEE_{n,H_n}(\Theta)$, provided that we have
$\Theta=[0,\infty)^4$ in \eqref{experiments}.

Next, we give a GARCH parametrization such that the randomly thinned
GARCH converges strongly to the MCOGARCH experiment $\widehat\EEE$ in
deficiency; therefore, we pick
$\theta=(h_0,\beta,\alpha,\lambda)\in\Theta$ and $n\in\NN$. If
$\alpha>0,$ then we set
%
\begin{eqnarray}\label{specialHn1}
h^{(0)}_{0,n}(\theta)&=&h_0\mathrm{e}^{-\alpha/n}+\frac
\beta\alpha(1 - \mathrm{e}^{-\alpha/n}),\qquad
\beta^{(0)}_n(\theta) = \frac\beta\alpha (1 - \mathrm{e}^{-\alpha
/n}) ,\nonumber\\ [-8pt]\\ [-8pt]
\alpha^{(0)}_n(\theta)&=&\mathrm{e}^{-\alpha/n},\qquad
\lambda^{(0)}_n(\theta) = \lambda\mathrm{e}^{-\alpha/n}\nonumber
\end{eqnarray}
and, otherwise, if $\alpha=0,$ then we set
%
\begin{equation}\label{specialHn2}
h^{(0)}_{0,n}(\theta)=h_0+\frac\beta
n,\qquad \beta^{(0)}_n(\theta) = \frac\beta n,\qquad  \alpha^{(0)}_n(\theta) = 1,\qquad \lambda
^{(0)}_n(\theta)=\lambda .
\end{equation}
Let $([0,\infty)^4,(H^{(0)}_n))$ be the corresponding parametrization
and $G^{(0)}_n$ the corresponding partial sum processes of GARCH in
\eqref{G1}.

Although the parametrization in \eqref{specialHn1} and \eqref{specialHn2}
is quite elaborate, we show that the corresponding GARCH experiment
converges to the experiment of MCOGARCH-type with no restrictions on
the limiting probability measure $Q$ assumed (see Section \ref{proofA}
for a~proof).
\begin{theorem}\label{th2}
Let \eqref{lre} be satisfied for
some $\gamma\in(0,\infty)$ and $p_n\in(0,1)$, $n\in\NN$.
If $Q_n$ tends to a probability measure $Q$ in total variation as
$n\to\infty,$ then $\EEE_{n,H_n^{(0)}}$ converges strongly to
$\widehat\EEE$ in deficiency as $n\to\infty$.
\end{theorem}

If $Q$ is absolutely continuous with respect to the Lebesgue measure,
then Theorem~\ref{th2} partially extends to other GARCH
parametrizations (see Section \ref{proofB} for a proof of the following
theorem).

\begin{theorem}\label{th2a} Let \eqref{lre}
be satisfied for some $\gamma\in(0,\infty)$ and $p_n\in(0,1)$,
$n\in\NN$. Suppose both that $Q_n$ tends to a probability measure $Q$
in total variation as $n\to\infty$ and that $Q\ll\ell$.

Let $\Theta\neq\varnothing$ with compact closure $\overline{\Theta}$ in
$(0,\infty)\times[0,\infty)^3$. For $n\in\NN$, let
$H_n=(h_{0,n},\beta_n,\alpha_n,\break\lambda_n)\dvtx\Theta\to[0,\infty)^4$
be a
GARCH parametrization and $G_n$ the corresponding GARCH model in~\eqref{G1}.

If there exist $n_0\in\NN$ and $C>0$ such that, for all $n\ge n_0$,
both
%
\begin{equation}\label{repunifbounded1}
\sup_{\theta=(h_0,\beta,\alpha,\lambda)
\in\Theta}\max \{ |h_{0,n}(\theta) - h_0 | ,
 |\lambda_{n}(\theta) - \lambda | \}\le\frac C n
\end{equation}
and
%
\begin{equation}\label{repunifbounded2}
\sup_{\theta=(h_0,\beta,\alpha,\lambda)
\in\Theta}\max \bigl\{ |n\beta_{n}(\theta) - \beta |,
\bigl|n\bigl(\alpha_n(\theta) - 1\bigr)+\alpha \bigr| \bigr\}\le C ,
\end{equation}
then
%
\begin{equation}\label{th2astronger}
\lim_{n\to\infty}\sup_{\theta
\in\Theta} \bigl\|
\LLL_{\theta}(G_n)-\LLL_{\theta}\bigl(G_n^{(0)}\bigr)  \bigr\|=0
\end{equation}
and $\EEE_{n,H_n}(\Theta)$ converges strongly to $\widehat
\EEE(\Theta)$ in deficiency as $n\to\infty$.
\end{theorem}
\begin{remark} \label{KVMMS}
(i) Let $Q=Q_n$ for all $n\in\NN$.
In Kallsen and Vesenmayer \cite{KV07} and Maller \textit{et al.}~\cite{MMS08},
the GARCH parametrizations $(\Theta,(H^{\mathrm{(KV)}}_n)_{n\in\NN})$ and
$(\Theta,(H^{\mathrm{(M)}}_n)_{n\in\NN})$ have been considered where, for
$\theta=(h_0,\beta,\alpha,\lambda)\in(0,\infty)^3\times[0,\infty)$,
using the obvious notation, $\Theta=(0,\infty)^3\times[0,\infty)$ and
%
\begin{eqnarray}\label{KVM}
h^{\mathrm{(KV)}}_{0,n}(\theta)&=&h^{\mathrm{(M)}}_{0,n}(\theta)=h_0 , \nonumber\\[-2pt]
\beta^{\mathrm{(KV)}}_n(\theta)&=&\beta^{\mathrm{(M)}}_n(\theta)=\frac\beta n
,\nonumber\\ [-9pt]\\ [-9pt]
\alpha^{\mathrm{(KV)}}_n(\theta)&=&\alpha^{\mathrm{(M)}}_n(\theta
)=\mathrm{e}^{-\alpha/n} ,\nonumber\\[-2pt]
\nonumber\lambda^{\mathrm{(KV)}}_n(\theta)&=&\lambda,\qquad
\lambda^{\mathrm{(M)}}_n(\theta) = \mathrm{e}^{-\alpha/n}\lambda .\nonumber
\end{eqnarray}
Kallsen and Vesenmayer \cite{KV07} have shown that
$(G_n[n\cdot],h_n[n\cdot])$, as defined in \eqref{G1} by
$H_n(\theta)=H^{\mathrm{(KV)}}_n(\theta)$, converge to COGARCH with parameter
$\theta$ in \eqref{COG} in law, with respect to the Skorokhod topology,
as $n\to\infty$, for all $\theta\in\Theta$.

Maller \textit{et al.} \cite{MMS08} have encountered a slightly different
scenario. For $\theta\in\Theta,$ they have embedded a sequence of GARCH
models into a given COGARCH and obtained the convergence with respect
to the same topology, now driven by a general L\'{e}vy process, even in
probability. If the driving process is a compound Poisson process with
rate $\gamma$ and jump size distribution $Q$, then their analysis
relates to a situation where the corresponding partial sums have the
same law as $(G_n[n\cdot],h_n[n\cdot])$ under the parametrization
$H^{\mathrm{(M)}}_n(\theta)$, $\theta\in\Theta$, $n\in\NN$.

In short, it follows from the analyses in \cite{KV07} that, for fixed
$\theta\in\Theta$, the partial sum
processes of GARCH with parametrization $H^{(KV)}_n(\theta)$ converge to
COGARCH in law, as
$n\to\infty$, with a similar result being true for the parametrization
$H^{(M)}_n(\theta)$ in~\cite{MMS08}. On the other hand,
both parametrizations fall into the framework of Theorem \ref{th2a}.
Hence, if the distribution of the innovations admits a Lebesgue
density, then the limiting experiment is given by MCOGARCH $\widehat
\EEE(\Theta)$ rather than COGARCH
$\EEE(\Theta)$.\vadjust{\goodbreak}

(ii)   In part (i), $Q=Q_n$ does not depend on $n$. Potential
applications where $Q_n$ depends on~$n$ arise in the
peak-over-threshold method in extreme value theory; see, for instance,
Embrechts \textit{et al.} \cite{EK97}, Resnick \cite{Re87} and
Falk \textit{et al.} \cite{FHR00}. Here, $Q_n$ equals the laws of
rescaled innovations, conditioned on the event that they exceed a given
threshold. Under reasonable assumptions, $Q_n$ converge weakly to a
generalized Pareto distribution~$Q$ as $n\to\infty$. Also, the
corresponding GARCH models converge in distribution in law to a COGARCH
driven by a compound Poisson process with jump distribution $Q$. In
this sense, COGARCH serves as a good approximation of GARCH \textit{in
law} if we are interested in the extreme parts of the innovations. On
the other hand, if $Q_n$ converges to $Q$, even in total variation
norm, then it follows from Theorem \ref{th2} that the corresponding
limiting experiment must be statistically equivalent to MCOGARCH.
\end{remark}

\subsection{COGARCH vs. MCOGARCH (incomplete
observations)}\label{secIncompletecomparison}
In this subsection, we
investigate whether the experiments induced by COGARCH and MCOGARCH are
of the same type. Here, we again assume that the volatility processes
are unobservable. Therefore, we recall \eqref{COG} and consider the
experiment
%
\begin{equation}\label{COGARCHexperiments}
\EEE =  \bigl(D_1,\DDD_1, (\LLL_{\theta}(G) )_{\theta\in
[0,\infty)^4} \bigr) .
\end{equation}
Note that both experiments $\EEE$ and $\widehat\EEE$ depend on the
intensity measure $\gamma\ell\otimes Q$ which enters~\eqref{COG} and
\eqref{hatEEE} via $N$. In this subsection, we include this dependence
in our notation by writing $\EEE_{\gamma,Q}$ and
$\widehat\EEE_{\gamma,Q}$ instead of $\EEE$ and $\widehat\EEE$,
respectively.

Let $f\dvtx\RR\to(0,\infty]$ be a strictly positive probability density
with respect to Lebesgue measure and set
%
\begin{equation}\label{defgQz}
g_{f, \zeta}(h):=h^\zeta  \int_\RR
f(hz)^\zeta f(z)^{1-\zeta}\,\mathrm{d}z,\qquad  h>0 , \zeta\in(0,1).
\end{equation}
By H\"{o}lder's inequality, $g_{f, \zeta}$ defines a function
$g_{f,\zeta}\dvtx(0,\infty)\to(0,1]$ with $g_{f,\zeta}(1)=1$. Note that
$g_{f, \zeta}$ satisfies both a scaling and a reflection property: for
all $0<\zeta<1$, $a,h>0$,
%
\begin{equation}\label{gfzscalingrefl}
g_{af(a\cdot), \zeta}(h)=g_{f, \zeta}(h),\qquad
g_{f, \zeta}(h)=g_{f, 1 - \zeta}(1/h).
\end{equation}
Next, we
investigate how COGARCH relates to MCOGARCH in deficiency (see
Section~\ref{proofth1} for a~proof).
\begin{theorem}\label{th1}
Let $\varnothing\neq\Theta\subseteq
(0,\infty)\times[0,\infty)^3$.
Assume that $Q$ admits a strictly positive Lebesgue density $f$ such
that for some $\zeta_0\in(0,1)$, $g_{f, \zeta_0}\dvtx(0,\infty)\to[0,1]$
is strictly increasing on $(0,1]$.

Let $(\gamma_n)_{n\in\NN}\subseteq(0,\infty)$ be a sequence such that
$\gamma=\lim_{n\to\infty} \gamma_n$ exists in $[0,\infty)$ and
$\gamma_n\neq\gamma$ for all $n\in\NN$.

If $\EEE_{\gamma_n,Q}(\Theta)$ is equivalent to
$\widehat\EEE_{\gamma_n,Q}(\Theta)$ for all $n\in\NN,$ then we
have:
\begin{longlist}
\item   if
$(h_{0,1},\beta,\alpha,\lambda),(h_{0,2},\beta,\alpha,\lambda)\in
\Theta$ and $\beta>0$, then $h_{0,1}= h_{0,2}$;
\item   if
$(h_{0},\beta_1,\alpha,\lambda),(h_{0},\beta_2,\alpha,\lambda)\in
\Theta$, then $\beta_{1}= \beta_{2}$;
\item   if
$(h_{0},\beta,\alpha_1,\lambda),(h_{0},\beta,\alpha_2,\lambda)\in
\Theta$ and $\beta=0$, then $\alpha_1=\alpha_{2}$;\vadjust{\goodbreak}
\item  if
$(h_{0},\beta,\alpha_1,\lambda),(h_{0},\beta,\alpha_2,\lambda)\in
\Theta$ and $\alpha_1=0$, then $\alpha_{2}=0$;
\item   if
$(h_{0},\beta,\alpha_1,\lambda),(h_{0},\beta,\alpha_2,\lambda)\in
\Theta$
and $\alpha_1<\alpha_2$, then $h_0>\beta/\alpha_{1}$.
\end{longlist}
\end{theorem}

Theorem 2.3 indicates that equivalence of MCOGARCH and COGARCH is
restricted to parameter sets that are of considerably lower dimension
and which have non-empty interiors. Hence, we do not have equivalence
in deficiency.

Observe that $\zeta\mapsto g_{f,\zeta}(h)$ occurs as the Hellinger
transformation of the scaling experiment
$(\RR,\BBB(\RR),\{\LLL(Z),\LLL(Z/h)\})$, where $Z$ is a random variable
with
Lebesgue density $f$. Next, we verify the monotonicity property of
$g_{f,\zeta}(h)$ in a number of examples.

\textit{Generalized symmetric gamma distribution.} Let $a,b,c>0$ and
$\Gamma$ be Euler's gamma function. Assume that $f\dvtx\RR\to(0,\infty]$
has the following form:
\[
f(z) =  \frac12 \frac{a^{c/b}}{\Gamma(c/b)} \mathrm{e}^{-a|z|^b}|z|^{c-1},\qquad  z\in\RR .
\]
This class of distributions covers important special cases such as the
normal distribution with zero mean and the Laplace distribution. It
follows straightforwardly that
\[
g_{f, \zeta}(h) =  \biggl(\frac{h^{b\zeta}}{h^b\zeta+(1-\zeta
)} \biggr)^{c/b},\qquad 0<\zeta<1,\  h>0 .
\]
Observe that $g_{f,\zeta}\dvtx(0,\infty)\to[0,1]$ is strictly increasing
on $(0,1]$ for all $0 < \zeta < 1$ and thus~$f$ satisfies the monotonicity assumption of Theorem \ref{th1}.

\textit{Centered Cauchy distribution.} Let $a>0$ and $f_a(z)= \frac
a{\uppi}\frac1{1+(az)^2}$ be the density of the centered $\operatorname{Cauchy}$
distribution $\operatorname{Cauchy}(0,a)$ with scaling parameter $a$. By the scaling
property in \eqref{gfzscalingrefl}, we have
$g_{f_a, 1/2}(h)=g_{f_1, 1/2}(h)$ for all $h>0$. By differentiating
this under the integral sign, we obtain
\[
\frac{\mathrm{d}}{\mathrm{d}h}g_{f_a, 1/2}(h)=\frac{1 - h^2}{2\uppi h^2} \int
_0^\infty
\frac{\sqrt{x}\,\mathrm{d}x}{(1+(1+h^2)x/h+x^2)^{3/2}}>0,\qquad
a>0 , 0<h\le1 .
\]
Consequently, the centered $\operatorname{Cauchy}$ distribution
satisfies the monotonicity assumption of Theorem \ref{th1}.

Next, we present a simulation-based approach to assess non-equivalence.
This approach can be used in cases not covered by Theorem \ref{th1} (or
when it is not clear whether the assumption of Theorem \ref{th1} is
satisfied). Recall that statistical equivalence of the experiments
$\EEE$ and $\widehat\EEE$ is implied (see \cite{St85}, Theorem 53.10)
when, for all finite subsets $\Theta\subseteq[0,\infty)^4$ and all
$\theta_0\in\Theta$, we have
%
\begin{equation}\label{COGMCOGdis}
\LLL_{\theta_0} \biggl( \biggl(\frac{\mathrm{d}\LLL_{\theta}(G)}{\mathrm{d}\LLL
_{\theta_0}(G)} \biggr)_{\theta\in\Theta} \biggr)
  =   \LLL_{\theta_0} \biggl( \biggl(\frac{\mathrm{d}\LLL_{\theta}(\widehat
G)}{\mathrm{d}\LLL_{\theta_0}(\widehat G)} \biggr)_{\theta\in\Theta} \biggr).
\end{equation}
We generated samples from these two distributions according to the
recursion (\ref{defmap1}) in the proof of Theorem \ref{th1} in Section
\ref{proofth1}. To this end, we first restricted the parameter space to
a set with two elements, $\theta_0$ and $\theta$. While fixing
$\theta_0$ to $(2,1,1,0.1)$, we have chosen eight vectors\vadjust{\goodbreak}
$\theta_{ij}$, $i=1,\ldots,4$, $j=1,2$, for the parameter vector
$\theta$, which differ from $\theta_0$ in only one component; see Table~\ref{TABchooseTHETA}.
%
\begin{table}
\tablewidth=210pt
\caption{Choices of $\theta_0$ and
$\theta=\theta_{ij}$ in equation (\protect\ref{COGMCOGdis})}\label{TABchooseTHETA}
\begin{tabular*}{210pt}{@{\extracolsep{\fill}}ld{2.1}lld{1.2}@{}}
\hline
$\theta_0$ & 2 & 1 & 1 & 0.1 \\ [3pt]
$\theta_{11}$ & 0.4 & 1 & 1 & 0.1 \\
$\theta_{12}$ & 10 & 1 & 1 & 0.1 \\ [3pt]
$\theta_{21}$ & 2 & 0.2 & 1 & 0.1 \\
$\theta_{22}$ & 2 & 5 & 1 & 0.1 \\ [3pt]
$\theta_{31}$ & 2 & 1 & 0.2 & 0.1 \\
$\theta_{32}$ & 2 & 1 & 5 & 0.1 \\ [3pt]
$\theta_{41}$ & 2 & 1 & 1 & 0.02 \\
$\theta_{42}$ & 2 & 1 & 1 & 0.5 \\
\hline
\end{tabular*}
\vspace*{-3pt}
\end{table}
Second, we checked the distributional equality (\ref{COGMCOGdis}) for
three different jump distributions: the~standard normal, the standard
$\operatorname{Cauchy}$ distribution $\operatorname{Cauchy}(0,1)$ (for comparison, note that
both of these are covered by Theorem \ref{th1}) and the normal mixture
distribution
\[
\tfrac{1}{2} N(-0.5,0.75) + \tfrac{1}{2} N(0.5,0.75),
\]
which has mean 0 and variance 1. The intensity $\gamma$ was always
fixed to $4$.

For each of the eight pairs $(\theta_0,\theta_{ij})$ and each of the
three jump distributions, we generated~$10^6$ samples of the two
distributions referring to the COGARCH and MCOGARCH in equation~(\ref{COGMCOGdis}).
The left-hand column of Table \ref{TABquantiles} reports the choice of
$\theta_{ij}$, whereas the other nine columns report, for each of the
three jump distributions, the 25\% quantile, median and 75\% quantile
of the distribution in equation (\ref{COGMCOGdis}).

%
\begin{table}
\tabcolsep=0pt
\caption{Estimated 25\% quantiles, medians and
75\% quantiles for the distributions in (\protect\ref{COGMCOGdis})}\label{TABquantiles}
\begin{tabular*}{\textwidth}{@{\extracolsep{\fill}}l@{}l@{}lll@{}l@{}lll@{}l@{}lll@{}}
\hline
& & \multicolumn{11}{@{}l@{}}{COGARCH}\\
& & \multicolumn{11}{@{}l@{}}{MCOGARCH}\\ [-7pt]
&&\multicolumn{11}{@{}l@{}}{\hrulefill}\\ [-3pt]
\multirow{1}{25pt}{Jumps} &\hspace*{10pt}& \multicolumn{3}{@{}l@{}}{$N(0,1)$} &\hspace*{10pt}& \multicolumn{3}{@{}l@{}}{$\operatorname{Cauchy}(0,1)$}
&\hspace*{10pt}& \multicolumn{3}{@{}l@{}}{$\mbox{Mixed }N$} \\ [-7pt]
&&\multicolumn{3}{@{}c@{}}{\hrulefill}&&\multicolumn{3}{@{}c@{}}{\hrulefill}&&\multicolumn{3}{@{}c@{}}{\hrulefill}\\
[-2pt]
Quantiles&& 25\%& Median& 75\% && 25\%& Median &75\% && 25\%&
Median &75\% \\
\hline
$\theta_{11}$ && 0.1081& 0.5560 &1.3888 && 0.5521 &0.7775 &1.1767 && 0.0909&
0.5329 &1.3918 \\
&& 0.1785 &0.6977 &1.3495 && 0.5884 &0.8226 &1.1811 && 0.1558& 0.6743& 1.3543 \\
[3pt]
$\theta_{12}$ && 0.1505 &0.3152& 0.6449 && 0.4173& 0.8127& 1.4573 && 0.1436&
0.3008 &0.6136 \\
&& 0.1637 &0.3377 &0.6768 && 0.4412 &0.8335 &1.4505 && 0.1575 &0.3264& 0.6559 \\
[3pt]
$\theta_{21}$ && 0.8326 &1.0168 &1.1711 && 0.9273& 0.9761 &1.0393 && 0.8307&
1.0201 &1.1766 \\
&& 0.7605 &1.0114 &1.2459 && 0.9051& 0.9566 &1.0539 && 0.7560 &1.0155 &1.2512 \\
[3pt]
$\theta_{22}$ && 0.4883 &0.7071 &1.0086 && 0.7765 &1.0229& 1.2130 && 0.4797&
0.6956 &1.0000 \\
&& 0.4201 &0.6077 &1.0000 && 0.7010& 1.0247 &1.2676 && 0.4100& 0.5988& 0.9798 \\
[3pt]
$\theta_{31}$ && 0.6928 &0.8543 &1.0621 && 0.8497& 1.0000 &1.1506 && 0.6863&
0.8476 &1.0530 \\
&& 0.6304 &0.7841 &1.0629 && 0.8029& 1.0000& 1.1881 && 0.6248& 0.7757& 1.0524 \\
[3pt]
$\theta_{32}$ && 0.0053 &0.1702 &1.1056 && 0.3853 &0.6449 &1.1172 && 0.0028&
0.1392 &1.0856 \\
&& 0.0010 &0.0590 &0.9129 && 0.3093& 0.5650& 1.1090 && 0.0005& 0.0437& 0.8703 \\
[3pt]
$\theta_{41}$ && 0.9864& 1.0104 &1.0735 && 0.8265& 1.0000& 1.0798 && 0.9863&
1.0114& 1.0762 \\
&& 0.9884 &1.0100& 1.0693 && 0.8357 &1.0000 &1.0779 && 0.9884 &1.0109 &1.0722 \\
[3pt]
$\theta_{42}$ && 0.6851& 0.8870 &1.0000 && 0.6217& 0.9328 &1.0418 && 0.6750&
0.8802& 1.0000 \\
&& 0.6963& 0.8942& 1.0000 && 0.6281 &0.9360 &1.0388 && 0.6865 &0.8874& 1.0000 \\
\hline
\end{tabular*}
\end{table}

Next, we applied the Wilcoxon rank sum test (also known as
Mann--Whitney test) to investigate the null hypothesis \textit{the median
of the likelihood ratio for the COGARCH experiment equals the median of
the likelihood ratio for the MCOGARCH experiment}. Table~\ref{TABWilcox} reports the values of the Wilcoxon test statistic $W$,
together with the corresponding $p$-values. For each jump distribution,
the first column corresponds to a sample size of $10^4$, the second to
$10^5$ and the third to a~sample size of $10^6$ per experiment.
Obviously, the $p$-values tend to 0 as the sample size increases. Based
on $10^6$ samples, the null hypothesis is most significantly rejected,
for all three jump distributions and all eight parameter vectors
$\theta_{ij}$. In other words, there is strong evidence that in the
case of incomplete observations, the randomly thinned GARCH and
COGARCH experiments are not statistically equivalent for these jump
distributions. This confirms our conjecture, that Theorem \ref{th1}
holds in a much more general formulation for quite arbitrary jump
distributions.

%
\begin{table}
\tabcolsep=1pt
\caption{Wilcoxon rank sum test: Values of Wilcoxon
test statistic $W$ and corresonding $p$-values}\label{TABWilcox}
{\fontsize{7.5}{10}\selectfont{\begin{tabular*}{\textwidth}{@{\extracolsep{\fill}}l@{}c@{}d{3.4}d{4.4}d{4.4}@{}c@{}d{3.4}d{4.4}d{4.4}@{}c@{}d{3.4}d{4.4}d{4.4}@{}}
\hline
& &\multicolumn{11}{l@{}}{W statistic}\\
& &\multicolumn{11}{l@{}}{$p$-value}\\ [-7pt]
&& \multicolumn{11}{l@{}}{\hrulefill}\\
\multicolumn{1}{@{}l@{}}{Jumps} &\hspace*{10pt}& \multicolumn{3}{@{}l@{}}{$N(0,1)$} &\hspace*{10pt}& \multicolumn{3}{@{}l@{}}{$\operatorname{Cauchy}(0,1)$}
&\hspace*{10pt}&
\multicolumn{3}{@{}l@{}}{$\mbox{Mixed }N$} \\ [-7pt]
&&\multicolumn{3}{@{}l@{}}{\hrulefill}&&\multicolumn{3}{@{}l@{}}{\hrulefill}&&\multicolumn{3}{@{}l@{}}{\hrulefill}\\
[-2pt]
\multicolumn{1}{@{}l@{}}{Sample size} && \multicolumn{1}{@{}l}{$10^4$} & \multicolumn{1}{l}{$10^5$} &
\multicolumn{1}{l@{}}{$10^6$}&
& \multicolumn{1}{@{}l}{$10^4$} & \multicolumn{1}{l}{$10^5$} & \multicolumn{1}{l@{}}{$10^6$}
&&
\multicolumn{1}{@{}l}{$10^4$}& \multicolumn{1}{l}{$10^5$}  & \multicolumn{1}{l@{}}{$10^6$}\\
\hline
$\theta_{11}$ && -7.10 & -24.12& -73.91 && -8.82& -25.46& -73.81 && -8.11&
-24.01 & -71.91 \\
&& 0.0000 & 0.0000 & 0.0000 && 0.0000&  0.0000 &0.0000 && 0.0000 &0.0000 &0.0000 \\
[3pt]
$\theta_{12}$ && -3.04 & -14.90&  -47.21 && -0.35 & -2.73&  -9.90 && -6.04 &
-15.37&
-48.40 \\
&& 0.0024 & 0.0000 & 0.0000 && 0.7245 & 0.0064& 0.0000 && 0.0000& 0.0000 & 0.0000 \\
[3pt]
$\theta_{21}$ && -1.56 & -3.20 & -12.52 && 8.71 & 31.90& 98.90 && -0.45 &-3.45&
-14.05 \\
&& 0.1189 & 0.0014& 0.0000 && 0.0000& 0.0000& 0.0000 && 0.6545& 0.0006& 0.0000 \\
[3pt]
$\theta_{22}$ && 12.10& 44.13& 136.09 && -1.92& -2.69& -8.28 && 14.17&
45.16&
141.15 \\
&& 0.0000 &0.0000& 0.0000 && 0.0546 &0.0070 &0.0000 && 0.0000 &0.0000& 0.0000 \\
[3pt]
$\theta_{31}$ && 12.38 &37.96 &116.09 && 1.76& 2.16 &10.30 && 12.07& 38.89&
119.95 \\
&& 0.0000& 0.0000 &0.0000 && 0.0788& 0.0311 &0.0000 && 0.0000 &0.0000& 0.0000 \\
[3pt]
$\theta_{32}$ && 11.34 &39.48& 126.96 && 11.63& 33.04 &100.34 && 13.59 &42.66&
131.75 \\
&& 0.0000 &0.0000& 0.0000 && 0.0000& 0.0000& 0.0000 && 0.0000& 0.0000 &0.0000 \\
[3pt]
$\theta_{41}$ && 0.83& 1.41& 5.79 && -1.52& -2.85& -4.98 && 2.83 &3.64& 5.35 \\
&& 0.4054 &0.1572 &0.0000 && 0.1280 &0.0044& 0.0000 && 0.0047 &0.0003& 0.0000 \\
[3pt]
$\theta_{42}$ && -2.96 &-3.71& -13.73 && -1.29 &-2.35 &-2.94 && -1.63 &-4.81&
-15.26 \\
&& 0.0031 &0.0002 &0.0000 && 0.1963 &0.0189& 0.0032 && 0.1041& 0.0000 &0.0000 \\
\hline
\end{tabular*}}}
\vspace*{3pt}
\end{table}
%

\vspace*{-2pt}\subsection{Complete observations}\vspace*{-2pt}\label{secComplete}
In the previous
subsections, we investigated both convergence and equivalence in
deficiency of a variety of GARCH-type experiments under the
assumption
that their volatility processes~$h_n$,\vadjust{\goodbreak} $h$ and $\hat h$ are
unobservable. In this subsection, we deal with the situation where the
corresponding volatility processes are observable in full. Of course,
this situation is mainly of theoretical interest and will primarily
help us to learn about the structural connections between GARCH and
COGARCH. However, we want to briefly mention some modern approaches by
which the unobservability of the volatility process can be dealt with
in practice. For example, there are several modern ways to estimate the
local volatility directly; see, for example, A{\"{i}}t-Sahalia, Mykland
and Zhang \cite{AMZh} and references therein or Jacod,
Kl\"{u}ppelberg and M\"{u}ller \cite{JaKlMu}, who use local
volatility estimates also in a COGARCH context, and many others. The
paper by Hubalek and Posedel~\cite{HubPos} contains another
very interesting idea. They use martingale estimating functions to
estimate the parameters in the Barndorff-Nielsen--Shephard model, which
is composed of a~stochastic differential equation (SDE) for the log
prices and another SDE for the variance. However, the martingale
estimating functions approach requires that both processes can be
observed. Hence, Hubalek and Posedel~\cite{HubPos} reinterpret the volatility
equation as an equation for some other observable measure of trading
intensity (such as trading volume or the number of trades), assuming
that the instantaneous variance process behaves (up to a
time-independent constant) exactly as the observable trading volume (or
the number of trades). As they show in their real-data example, this
approach leads to quite satisfying results. The same idea could be
used, of course, for the COGARCH model, to bypass problems with the
unobservability of the volatility process in practice.

Returning to theoretical matters, we now consider the following
GARCH-type experiments in continuous time with fully observed
volatilities, denoted by
\[
\EEE_h= \bigl(D_2,\DDD_2, (\LLL_{\theta}(G,h) )_{\theta\in
[0,\infty)^4} \bigr),\qquad
\widehat\EEE_h= \bigl(D_2,\DDD_2, (\LLL_{\theta}(\widehat
G,\hat
h) )_{\theta\in[0,\infty)^4} \bigr) ,
\]
where $\hat h$ is defined by the specification in \eqref{hatEEE}
and \eqref{timechange}. Similarly to
Sections \ref{secGARCH} and~\ref{secIncomplete}, where we dealt with
continuous time, both experiments $\EEE_h$ and $\widehat\EEE_h$
depend on $Q$ and $\gamma>0$ as well. In this subsection, we will
suppress this dependence in our notation.

We need to specify a set $\Theta_{e}\subseteq[0,\infty)^4$ of {\em
exceptional points} in the parameter space $[0,\infty)^4$ as follows:
%
\begin{equation}\label{defThetae}
\Theta_{e} = \{\theta=(h_0,\beta,\alpha,\lambda)\in[0,\infty
)^4\dvtx h_0\alpha=\beta\} .
\end{equation}
Observe that $\Theta_e$ is closely connected to the fixed point of the
affine differential equation $h'(t)=\beta-\alpha h(t)$. Indeed, if
$\theta=(h_0,\beta,\alpha,\lambda)\in\Theta_{e}$, then we have
$h(t)=\hat h(t)\equiv h_0$ for all $t\in[0,T)$, where $T$ is the
first jump of (M)COGARCH. It is impossible to recover the parameters
$\beta,\alpha,\lambda$ in full within the time horizon $[0,T)$.
Otherwise, if $h_0$ is not the fixed point of this differential equation,
then it is always possible to recover parts of $\theta$ by taking
appropriate derivatives. In the next proposition, we formalize this
idea and show that both $\EEE_h$ and $\widehat\EEE_h$ are equivalent
to a simple reference experiment (see Section \ref{proofprop1} for a~proof).

\begin{proposition}\label{prop1}
If $Q(\{0\})=0,$ then
both $\EEE_h$ and $\widehat\EEE_h$ are equivalent to
$\FFF= ([0,\infty]^4,\BBB([0,\break\infty]^4),(Q_{\theta})_{\theta
\in[0,\infty)^4} )$
where, for $\theta=(h_0,\beta,\alpha,\lambda)\in[0,\infty)^4$,
$\gamma>0$, we set\vspace*{-2pt}
%
\begin{equation}\label{defhatQthetan}
Q_{\theta} =  \cases{
\mathrm{e}^{-\gamma}
\varepsilon_{(h_0,\beta,\alpha,\infty)} + (1 - \mathrm{e}^{-\gamma})
\varepsilon_{\theta},&\quad $\theta\notin\Theta_{e}$,\cr
\mathrm{e}^{-\gamma}
\varepsilon_{(h_{0},\infty,\infty,\infty)} + (1 - \mathrm{e}^{-\gamma})
\varepsilon_{\theta} ,&\quad $\theta\in\Theta_{e} , h_0 > 0,\lambda > 0$,\cr
\mathrm{e}^{-\gamma}
\varepsilon_{(h_{0},\infty,\infty,\infty)} + (1 - \mathrm{e}^{-\gamma})
\varepsilon_{(h_{0},\infty,\infty,0)} ,&\quad $\theta\in\Theta_{e},
h_0 > 0 , \lambda = 0$,\cr
\varepsilon_{(0,\infty,\infty,\infty)} ,&\quad $\theta\in\Theta_{e}, h_0 = 0$,
}\vspace*{-2pt}
\end{equation}
and $\Theta_e$ is the set defined as in \eqref{defThetae}.\vspace*{-2pt}
\end{proposition}
\begin{remark}\label{whyQonoto}
In the situation of Proposition \ref
{prop1}, we
require $Q$ to satisfy $Q(\{0\})=0$. Indeed, if $Q=\varepsilon_0$, then
it is easy to see that both $\EEE_h$ and $\widehat\EEE_h$ are
equivalent to $\FFF$, where we formally set $\gamma=0$ in
\eqref{defhatQthetan}. Otherwise, if $Q(\{0\})\in[0,1)$, then we may
\vspace{1pt}
adjust the intensity measures of the driving Poisson measure
accordingly, to see that both $\EEE_h$ and $\widehat\EEE_h$ are
equivalent to $\FFF$, but with $\gamma$ replaced by
$\gamma Q(\RR\backslash\{0\})$ in the definition of $Q_\theta$.
Analogously, one can adjust the discrete-time experiments that we
consider in Proposition \ref{prop2}. We leave the details to the
reader.\vspace*{-2pt}
\end{remark}

Next, we investigate the discrete-time experiments. Note that the
initial value of $h$ is observable in continuous time. As a result, it
is always possible to recover the parameter~$h_0$ in full. To account
for this phenomenon in discrete time, we introduce the following
sequence of experiments, $\EEE_{h,n,H_n}$, indexed by $n\in\NN$, where\vspace*{-2pt}
%
\begin{equation}\label{experimentsh2}
\EEE_{h,n,H_n} =  \bigl([\RR^{n + 1}]^2,\BBB([\RR^{n + 1}]^2),
 (\LLL_{\theta} (G_n,h_n ) )_{\theta\in[0,\infty
)^4} \bigr),\qquad n\in\NN .\vspace*{-2pt}
\end{equation}
Here, $([0,\infty)^4,(H_n))$ is a parametrization of the full parameter
space $[0,\infty)^4$; both $G_n=(G_{n,k})_{0\le k\le n}$
and $h_n=(h_{n,k})_{0\le k\le n}$ are defined by \eqref{G1} via
$H_n(\theta)=(h_{0,n}(\theta),\beta_n(\theta),\alpha_n(\theta
),\break\lambda_n(\theta))$
for $n\in\NN$ and $\theta\in[0,\infty)^4$ (by a slight abuse of
previous notation). We are now in a~position to state a discrete-time
analog of Proposition \ref{prop1} (see Section \ref{proofprop2} for a
proof).

\begin{proposition}\label{prop2}Suppose that \eqref{lre} is
satisfied for some $\gamma\in(0,\infty)$ and $p_n\in(0,1)$, $n\in
\NN$.
Let $([0,\infty)^4,H_n)_{n\in\NN}$ be the parametrization in
\eqref{specialHn1} and \eqref{specialHn2}. Also, let
$([0,\infty)^4,\allowbreak H_n^{\mathrm{(KV)}})_{n\in\NN}$ and
$([0,\infty)^4,H_n^{\mathrm{(M)}})_{n\in\NN}$ be the parametrizations in
\eqref{KVM}, respectively.\looseness=-1

If $Q(\{0\})=Q_n(\{0\})=0$ for all $n\in\NN$, then the following
assertions hold as $n\to\infty$, both in deficiency:
\begin{longlist}
\item   $\EEE_{h,n,H_n}$ converges strongly to
$\FFF$;
\item   both $\EEE_{h,n,H_n^{\mathrm{(KV)}}}$ and $\EEE_{h,n,H_n^{\mathrm{(M)}}}$ are
asymptotically equivalent to\vspace*{-2pt}
\[
\FFF_{n}= \bigl([0,\infty]^4,\BBB([0,\infty]^4),
(Q_{\theta,n})_{\theta\in[0,\infty)^4} \bigr),\vspace*{-2pt}
\]
where, for $n\in\NN$ and
$\theta=(h_0,\beta,\alpha,\lambda)\in[0,\infty)^4$, we define
$Q_{\theta,n}$ as $Q_{\theta}$ in \eqref{defhatQthetan}, but with~$\Theta_e$ replaced by\vspace*{-2pt}
\[
\Theta_{e,n} = \{\theta=(h_0,\beta,\alpha,\lambda)\in[0,\infty
)^4\dvtx
h_0 n (1 - \mathrm{e}^{-\alpha/n})=\beta\} .\vspace*{-2pt}
\]
\end{longlist}
\end{proposition}

Finally, we are concerned with the relationships between the
experiments $\FFF$ and $\FFF_n$, $n\in\NN\cup\{\infty\}$ (see
Section \ref{proofprop3} for a proof).\vadjust{\goodbreak}
\begin{proposition}\label{prop3} Let $\gamma>0$ and $\varnothing\neq
\Theta\subseteq[0,\infty)^4$. Let $\FFF, \FFF_n$, $n\in\NN$,
be the experiments in Propositions \ref{prop1} and \ref{prop2},
respectively.

Let $\widehat\FFF= ([0,\infty]^4,\BBB([0,\infty]^4),(\widehat
Q_\theta)_{\theta\in[0,\infty)^4} )$ be the experiment where, for
$\theta=(h_0,\beta,\alpha,\lambda)\in[0,\infty)^4$, we define
$\widehat
Q_{\theta}$ as $Q_{\theta}$ in \eqref{defhatQthetan}, but with
$\Theta_e$ replaced by
\[
\widehat\Theta_{e} = \{0\}^2\times(0,\infty)\times[0,\infty
)\cup[0,\infty)\times\{0\}^2\times[0,\infty) .
\]
The following assertions then hold:
\begin{longlist}
\item
$\delta(\widehat\FFF(\Theta),\FFF(\Theta))=\delta(\widehat\FFF
(\Theta),\FFF_n(\Theta))=0$
for all $n\in\NN$;
\item $\delta(\FFF(\Theta),\widehat\FFF(\Theta))=0$ if and
only if, for all $h_0>0$,
%
\begin{eqnarray}\label{hsgleichbetaalphagleich}
 &&\{(\beta,\alpha,\lambda)\in[0,\infty)^3\dvtx
(h_0,\beta,\alpha,\lambda)\in\Theta\cap\Theta_e\cap
\widehat\Theta_e^C \}\neq\varnothing\nonumber\\ [-10pt]\\ [-10pt]
&&\quad\Rightarrow\quad
\# \{(\beta ,\alpha)\in[0,\infty)^2 \dvtx \exists\lambda\ge
0\ (h_0,\beta,\alpha,\lambda)\in\Theta_e\cap\Theta \}=1;\nonumber
\end{eqnarray}
\item   $\lim_{n\to\infty}
\delta(\FFF_n(\Theta),\widehat\FFF(\Theta))=0$ if and only if there
exists some $n_0$ such that for all $n\ge n_0$ and $h_0>0$,
\eqref{hsgleichbetaalphagleich} holds, but with $\Theta_{e}$ replaced
by $\Theta_{e,n}$. In particular, $\FFF_n$ converges weakly to
$\widehat\FFF$ as $n\to\infty$ in deficiency.
\end{longlist}
\end{proposition}

Let us rephrase our results in terms of the GARCH experiments, with the
volatility processes fully observed in both continuous and discrete
time. In contrast to the situation in Theorem \ref{th1}, it follows
from Proposition \ref{prop1} that the continuous-time experiments
induced by (M)COGARCH are mutually equivalent in deficiency. Depending
on the parametrization, (M)COGARCH also occurs as the limit in
deficiency of discrete-time GARCH; in particular, this is the case for
the parametrization in Proposition \ref{prop2}. In contrast to
Theorem \ref{th1}, for a large class of parameter sets~$\Theta$, all of
these discrete-time experiments, that is,
$\EEE_{h,n,H^{(0)}_n}(\Theta)$, $\EEE_{h,n,H^{\mathrm{(KV)}}_n}(\Theta)$,
$\EEE_{h,n,H_n^{\mathrm{(M)}}}(\Theta)$, are asymptotically equivalent to
(M)COGARCH $\EEE_h(\Theta)$ and $\widehat\EEE_h(\Theta)$, in
deficiency, as $n\to\infty$. This happens, for instance, if
$\Theta\subseteq[0,\infty)^4$ does not contain an open neighborhood of
$\Theta_{e}$. Since the set~$\Theta_e$ is of lower dimension than
$[0,\infty)^4$ it is thus justified to say that the randomly thinned
GARCH is \textit{generically} equivalent to COGARCH in deficiency as
$n\to\infty$.

\vspace*{-3pt}\section{Conclusion}\vspace*{-3pt}\label{summary}
In Le Cam's framework, Wang \cite{Wa02} and Brown \textit{et al.} \cite{BWZ03}
investigated GARCH and Nelson's diffusion limit. These authors dealt
with aggregated Gaussian innovations. For a~suitable parametrization,
Maller \textit{et al.} \cite{MMS08} and Kallsen and Vesenmayer \cite{KV07}
showed that the GARCH model converges to the COGARCH model in
probability and in distribution, respectively, when the innovations are
randomly thinned. These papers dealt with a general L\'{e}vy process as
the driving process of the COGARCH. In this paper, we have studied
an
important special case in Le Cam's framework of statistical
experiments, namely, we have assumed that the driving process of
CO\-GARCH is a compound Poisson process.
GARCH then converges generically to COGARCH, even in deficiency,
provided that the volatility processes are observed. Hence, from a
theoretical point of view, COGARCH can indeed be considered as a
continuous-time equivalent to GARCH.\vadjust{\goodbreak} Otherwise, when the observations
are incomplete, GARCH still has a limiting experiment, which we call
MCOGARCH, but this will usually not be equivalent to COGARCH in
deficiency. Nevertheless, this limiting experiment is, from a
statistical point of view, quite similar to COGARCH since the only
difference is the exact localization of the jump times. For COGARCH,
the jump times can be more random than for MCOGARCH, but practitioners
may see this as an additional advantage of COGARCH.

It would be interesting to extend the analysis to more general L\'{e}vy
processes, rather than Brownian motion and compound Poisson processes.
However, this would first require substantial investigations of the
approximation and randomizations of L\'{e}vy processes themselves and,
therefore, seems out of reach at the present stage of research.

%
\vspace*{-3.5pt}\section{\texorpdfstring{Proof of Theorem \protect\ref{th2}}
{Proof of Theorem 2.1}}\vspace*{-3.5pt}\label{proofA}
For the reader's convenience, we first provide a brief roadmap for the
proof of Theorem~\ref{th2}. The proof is split into two parts, which
appear in Sections \ref{proofth2} and \ref{SUBSECproofTh2PartII},
respectively. The second part uses a lemma which we formulate and
prove in Section \ref{SUBSEC4.2}. To prove that $\EEE_{n,H_n^{(0)}}
\to
\hat{\EEE}$ in deficiency, we will introduce intermediate experiments
$\EEE^{\star}_{1,n}$ and $\EEE^{\star}_{2,n}$. The first of\vspace*{1pt} these two
experiments corresponds to a deterministic time grid, the latter to a
randomized time grid. First, we will show that $\EEE_{n,H_n^{(0)}}$ is
equivalent to $\EEE^{\star}_{1,n}$ in deficiency and then, using Lemma
\ref{thinning} from Section~\ref{SUBSEC4.2}, that $\EEE^{\star}_{2,n}$ converges
strongly to~$\widehat\EEE$. Finally, we will prove that
$\EEE^{\star}_{1,n}$ and $\EEE^{\star}_{2,n}$ are equivalent.
%
\vspace*{-3.5pt}\subsection{\texorpdfstring{Proof of Theorem \protect\ref{th2} (part I)}
{Proof of Theorem 2.1 (part I)}}\vspace*{-3.5pt}\label{proofth2}
For $n\in\NN$, define a point measure $N_{1,n}$ on $[0,1]\times\RR$ by
%
\begin{equation}\label{defNn}
N_{1,n}=\sum_{k=1}^n 1_{Z_{n,k}\neq
0}  \varepsilon_{(k/n , Z_{n,k})},\qquad  n\in\NN .
\end{equation}
Using $N_{1,n}$, we pass from discrete to continuous time. For
$n\in\NN$, define
\[
\EEE^{\star}_{1,n}=\bigl(D_1,\DDD_1,(\LLL_\theta(G_{1,n}))_{\theta
\in[0,\infty)^4}\bigr),
\]
where, for all $0\le t\le1$, $n\in\NN$ and
$\theta=(h_0,\beta,\alpha,\lambda)\in[0,\infty)^4$,
$(G_{1,n},h_{1,n})$ is the unique pathwise solution of
the following system of integral equations ($t\in[0,1]$):
%
\begin{eqnarray}\label{G1p}
  G_{1,n}(t) & = & \int_{[0,t]\times\RR} h^{1/2}_{1,n}(s-) z
N_{1,n}(\mathrm{d}s,\mathrm{d}z) ,\nonumber\\ [-10pt]\\ [-10pt]
 h_{1,n}(t) & = & h_0 + \int_{[0,t]} \beta- \alpha
h_{1,n}(s-)\,\mathrm{d}s
+ \lambda\int_{[0,t]\times\RR} h_{1,n}(s-)z^2 N_{1,n}(\mathrm{d}s,\mathrm{d}z).\nonumber
\end{eqnarray}
Fix $\theta=(h_0,\beta,\alpha,\lambda)\in[0,\infty)^4$ with
$\alpha\neq
0$. By solving the linear ODE for $h_{1,n}$ in \eqref{G1p}, observe
that
%
\begin{equation}\label{ode}
h_{1,n}(t) = \frac\beta{\alpha}
 \bigl[1 - \mathrm{e}^{-\alpha[t-(k - 1)/n]} \bigr]+
\mathrm{e}^{-\alpha[t-(k - 1)/n]} h_{1,n} \biggl(\frac{k-1} n
\biggr)\vadjust{\goodbreak}
\end{equation}
for $(k - 1)/n\le t<k/n$, $1\le k\le n$ and $n\in\NN$. It thus
follows from \eqref{specialHn1} and \eqref{ode} that for all $n\in
\NN$,
\begin{eqnarray*}\label{inth2h2n}
h_{1,n}(1/n-)&=&h_0\mathrm{e}^{-\alpha/n}+\frac{\beta}{\alpha}[1 -
\mathrm{e}^{-\alpha/n}] = h_{0,n}(\theta),\\
h_{1,n} \biggl(\frac k n- \biggr)&=&\beta_n(\theta)+
h_{1,n} \biggl(\frac{k - 1}
n- \biggr)  [\alpha_n(\theta)+\lambda_n(\theta)
Z_{n,k-1}^2 ],\qquad 2\le
k\le n.
\end{eqnarray*}
In view of \eqref{G1} and the identities in the last display, we
thus have
$h_{n}(k)=h_{1,n} (((k + 1)/n)- )$ for all $n\in\NN$,
$0 \le  k\le n - 1$ and
$\theta=(h_0,\beta,\alpha,\lambda)\in[0,\infty)^4$ with $\alpha
>0$. A
similar argument is applicable to \eqref{specialHn2} and
$\theta=(h_0,\beta,0,\lambda)\in[0,\infty)^4$. It thus follows from
\eqref{G1} and \eqref{G1p} that
\[
\LLL_\theta \bigl(\bigl(G_{1,n}(k/n)\bigr)_{1\le k\le n} \bigr) =
\LLL_\theta ((G_{n}(k))_{1\le k\le n} ),\qquad
n\in\NN , \theta\in[0,\infty)^4 .
\]
Note that $G_{1,n}$ is constant
on $[(k  -  1)/n,k/n)$, $1 \le  k \le  n$ and $n\in\NN$. Hence,
$\EEE_{n,H_n^{(0)}}$ is equivalent to $\EEE^{\star}_{1,n}$ in
deficiency for all $n\in\NN$ by \eqref{defMarkov} and the monotonicity
theorem for\vspace*{1pt} Markov kernels (see~\cite{Re93}, Lemma 1.4.2(i)).

Next, we randomize the deterministic time grid. Therefore, let
$(U_k)_{k\in\NN}$ be an i.i.d. sequence of random variables independent
of the vector $Z_n$, where $U_k$ is uniformly distributed on $[0,1]$.
Set
%
\begin{equation}\label{defVnk}
V_{n,k} = \bigl((k - 1) + U_k\bigr)/n,\qquad 1\le k\le n,
\end{equation}
and define a point process $ N_{2,n}$ by
%
\begin{equation}\label{deftildeN}
N_{2,n}=\sum_{k=1}^n
1_{Z_{n,k}\neq0} \varepsilon_{(V_{n,k},Z_{n,k})},\qquad
n\in\NN .
\end{equation}
Let $T$ be as in \eqref{timechange}. For
$n\in\NN$, let
$\EEE^{\star}_{2,n}=(D_1,\DDD_1,(\LLL_{\theta}(G_{2,n}))_{\theta
\in[0,\infty)^4})$,
where for all\break $0\le t\le1$, $n\in\NN$ and
$\theta=(h_0,\beta,\alpha,\lambda)\in[0,\infty)^4$,
$(G_{2,n},h_{2,n})$ is the pathwise unique solution of the
following system of integral equations:
%
\begin{eqnarray}\label{EEE3}
 \hspace*{-15pt}G_{2,n}(t)& = &\int_{[0,t]\times\RR}
h^{1/2}_{2,n}(s-)z N_{2,n}(\mathrm{d}s,\mathrm{d}z),\nonumber\\ [-8pt]\\ [-8pt]
\hspace*{-15pt}h_{2,n}(t) & = & h_0+ \int_{[0,t]}\beta-\alpha
h_{2,n}(s-)\,\mathrm{d}T_{N_{2,n}} (s) + \lambda\int_{[0,t]\times\RR}
h_{2,n}(s-)z^2 N_{2,n}(\mathrm{d}s,\mathrm{d}z).\nonumber
\end{eqnarray}
To proceed with the proof of Theorem \ref{th2}, we need the following
lemma.

%
\subsection{\texorpdfstring{Statement and proof of Lemma \protect\ref{thinning}}
{Statement and proof of Lemma 4.1}}\label{SUBSEC4.2}

\begin{lemma}\label{thinning} Let $N$ be a Poisson measure with
intensity measure
$\gamma\ell\otimes Q$ and $N_{2,n}$ as in \eqref{deftildeN}. Suppose
that \eqref{lre} holds. If $Q_n$ tends to $Q$ in total variation as
$n\to\infty$, then $\lim_{n\to\infty}\|\LLL( N_{2,n})-\LLL(N)\|=0$.
\end{lemma}

\begin{pf}
Suppose that \eqref{lre} is satisfied for $n\in\NN$, $p_n\in(0,1)$ and
$\gamma\in(0,\infty)$. Let $B_{n,1}, \dots,B_{n,n}$ be independent
Bernoulli variables with parameter $p_n$. Suppose that
$(U_k,\zeta_{n,k})_{k\in\NN}$ is an i.i.d. sequence of random vectors
with independent components, where~$U_k$ is uniformly distributed on
$(0,1)$ and $\LLL(\zeta_{n,k})=Q_n$. Suppose that $B_{n,1}, \dots,
B_{n,n}$ and $(U_k,\zeta_{n,k})_{k\in\NN}$ are independent. Observe
that
\[
\LLL(N_{2,n}) = \LLL \Biggl(\sum_{k=1}^n B_{n,k} \varepsilon
_{(V_{n,k},\zeta_{n,k})} \Biggr)
\]
with $V_{n,k}=(k - 1+ U_k)/n$ for all $n\in\NN$ and\vspace*{1pt} $1\le k\le n$.

Let $\widehat N_n$ be a Poisson measure on $[0,1]\otimes\RR$ with
intensity measure $np_n\ell\otimes Q_n$ and define
\[
\widehat N_{n,k}(B)=\widehat N_n \biggl(B\cap \biggl( \biggl(\frac{k -
1}n,\frac kn \biggr]\times\RR \biggr) \biggr),\qquad  B\in\BBB
([0,1]\times\RR) .
\]
$\widehat N_{n,1},\dots,\widehat N_{n,n}$ are then independent Poisson
point processes, where for all $n\in\NN$, \mbox{$1\le k\le n$}, $\widehat
N_{n,k}$ has intensity measure
\[
n p_n[\ell\otimes Q_n]  \biggl(B\cap \biggl( \biggl(\frac{k -
1}n,\frac kn \biggr]\times\RR \biggr) \biggr),\qquad
  B\in\BBB([0,1]\times\RR) .
\]
By the monotonicity theorem for
Markov kernels (see \cite{Re93}, Lemma 1.4.2(i)), observe that for all
$n\in\NN$,
%
\begin{equation}\label{thin1}
\|\LLL(N_{2,n})-\LLL(\widehat N_n)\|
\le
 \Biggl\|\bigotimes_{k=1}^n\LLL \bigl(B_{n,k}
\varepsilon_{(V_{n,k},\zeta_{n,k})} \bigr)-\bigotimes_{k=1}^n\LLL
 (\widehat
N_{n,k}  ) \Biggr\| .
\end{equation}
Denote the Hellinger distance
between two probability measures $P_1$ and $P_2$ by $H(P_1,P_2)$. This
gives us the following upper bound (see \cite{Re93}, Section~1.3,
equation (1.23), and Section~1.3, equation (1.25)):
%
\begin{eqnarray}\label{thin2}
&&\Biggl\|\bigotimes_{k=1}^n\LLL \bigl(B_{n,k}
\varepsilon_{(V_{n,k},\zeta_{n,k})} \bigr)-\bigotimes_{k=1}^n\LLL
 (\widehat N_{n,k}) \Biggr\|\nonumber\\
&&\quad\le H \Biggl(\bigotimes_{k=1}^n\LLL \bigl(B_{n,k}
\varepsilon_{(V_{n,k},\zeta_{n,k})} \bigr),
\bigotimes_{k=1}^n\LLL (\widehat N_{n,k}) \Biggr)\\
&&\quad\le \Biggl(\sum_{k=1}^n H^2 \bigl(\LLL \bigl(B_{n,k}
\varepsilon_{(V_{n,k},\zeta_{n,k})} \bigr),\LLL (\widehat
N_{n,k} ) \bigr) \Biggr)^{1/2}.\nonumber
\end{eqnarray}
Fix $n\in\NN$ and $1\le k\le n$. Let
$(V_{n,k,l},\zeta_{n,k,l})_{l\in\NN}$ be an i.i.d. sequence of random
vectors with $\LLL(V_{n,k,l},\zeta_{n,k,l})=\LLL(V_{n,k})\otimes Q_n$,
$l\in\NN$. Suppose that $(V_{n,k,l})$ is independent of $B_{n,k}$ and
$\tau_{n,k}$, where $\tau_{n,k}$ is a Poisson\vadjust{\goodbreak} variable with parameter
$p_n$. We then have the following identities:
\[
\LLL \bigl(B_{n,k}
\varepsilon_{(V_{n,k},\zeta_{n,k})} \bigr)   =
\LLL \Biggl(\sum_{l=1}^{B_{n,k}}
\varepsilon_{(V_{n,k,l},\zeta_{n,k,l})} \Biggr),\qquad  \LLL
(\widehat N_{n,k})   =   \LLL \Biggl(\sum_{l=1}^{\tau_{n,k}}
\varepsilon_{(V_{n,k,l},\zeta_{n,k,l})} \Biggr) .
\]
By \cite{Re93}, Lemma 1.4.2(ii), for $n\in\NN$ and $1\le k\le n$, we
must have
%
\begin{equation}\label{thin5}
H \bigl(\LLL \bigl(B_{n,k}
\varepsilon_{(V_{n,k},Z_{n,k})} \bigr),\LLL (\widehat
N_{n,k} ) \bigr)  \le H(\LLL(B_{n,k}),\LLL(\tau_{n,k})) .
\end{equation}
As $H(\LLL(B_{n,k}),\LLL(\tau_{n,k}))\le3^{1/2}p_n$ (see
\cite{Re93}, Theorem 1.3.1(ii)), it follows from
\eqref{thin1}--\eqref{thin5} and~\eqref{lre} that
%
\begin{equation}\label{thin6}
\limsup_{n\to\infty}\|\LLL(\tilde N_n)-\LLL(\hat
N_n)\| \le \limsup_{n\to\infty}(3np_n^2)^{1/2}=0 .
\end{equation}
In view of a well-known upper bound of the laws of Poisson point
measures in terms of the corresponding intensity measures (see
\cite{Re93}, Section 3.2, equation (3.8)), it follows from \eqref{lre}
and $\|\ell\otimes Q\|=1$ that
\begin{eqnarray*}\label{thin7}
\|\LLL(\widehat N_n)-\LLL(N)\| & \le&
3\|\gamma\ell\otimes Q-np_n\ell\otimes Q_n\|\\
& \le& 3 |np_n-\gamma| + 3 np_n \|Q-Q_n\| \to0\qquad \mbox{as }
n\to\infty .
\end{eqnarray*}
By means of \eqref{thin6} and \eqref{thin7}, this completes the proof
of the lemma.
\end{pf}

%
\subsection{\texorpdfstring{Proof of Theorem \protect\ref{th2} (part II)}{Proof of Theorem 2.1 (part II)}}\label{SUBSECproofTh2PartII} 
Let $N$ be a Poisson measure with intensity measure $\gamma\ell
\otimes
Q$. It follows from \eqref{hatEEE} and~\eqref{EEE3} that there exists a
family of deterministic Markov kernels
$K_\theta\dvtx\MM_1\times\DDD_1\to[0,1]$, indexed by
$\theta\in[0,\infty)^4$, such that both
$\LLL_\theta(G_{2,n})=K_\theta\LLL(N_{2,n})$ and $\LLL_\theta
(\widehat
G)=K_\theta\LLL(N)$ for all $n\in\NN$ and $\theta\in\Theta$.
Since we
have assumed \eqref{lre}, the assertion of Lemma \ref{thinning} holds
and we thus get from \eqref{deflecamsmallerTV} and the monotonicity
theorem for Markov kernels (see \cite{Re93}, Lemma 1.4.2(i)) that as
$n\to\infty$,
\[
\Delta(\widehat\EEE,\EEE^{\star}_{2,n})\le\sup_{\theta\in
[0,\infty)^4}\|
\LLL_\theta(\widehat G)-\LLL_\theta(G_{2,n})\|\le
\|\LLL(N)-\LLL(N_{2,n})\|\to0 .
\]
Consequently, $\EEE^{\star}_{2,n}$
converges (strongly) to $\widehat\EEE$ in deficiency as $n\to\infty$.
Recall that $\EEE_{n,H_n^{(0)}}$ is equivalent to $\EEE_{1,n}^{\star}$
in deficiency for all $n\in\NN$. To complete the proof of the theorem,
it thus suffices to show that $\EEE^{\star}_{1,n}$ is equivalent to
$\EEE^{\star}_{2,n}$.

Therefore, let $\MM_0$ denote the space of all non-negative point
measures on $[0,1]$ with finite support. We equip this space with the
$\sigma$-algebra $\MMM_0$ generated by point evaluations (see
Reiss~\cite{Re93}, pages 5--6). Let $\MM_{0,1}\subseteq\MM_0$
be the subset of point measures $\sigma\in\MM_0$ such that there exist
$m\in\NN$ and $0 = t_0 < t_1 < \cdots < t_{m} <  1$ with
$\sigma=\sum_{k=1}^m\varepsilon_{t_k}$. For $\sigma\in\MM_0$, we define
mappings $T_{1,\sigma},T_{2,\sigma}\dvtx[0,1]\to[0,\infty)$ and
$T_{3,\sigma},T_{4,\sigma}\dvtx[0,1]\times\RR\to[0,\infty)\times\RR
$ as
follows: if $\sigma\in\MM_{0}\backslash\MM_{0,1}$, then for all
$t\in[0,1]$ and $x\in\RR$, we set $T_{1,\sigma}(t)=T_{2,\sigma}(t)=t$
and $T_{3,\sigma}(t,x)=T_{4,\sigma}(t,x)=(t,x)$. Otherwise, if
$\sigma\in\MM_{0,1}$, then there exist $m\in\NN$ and
$0 = t_0 < t_1 < \cdots < t_{m} < 1$ with
$\sigma=\sum_{k=1}^m\varepsilon_{t_k}$ and we set
\begin{eqnarray*}
T_{1,\sigma}(t)&=& \frac{t-t_k}{m(t_k-t_{k - 1})}+\frac km ,\qquad
t\in[t_{k-1},t_k),\ 1\le k\le m ,\\
\nonumber T_{1,\sigma}(t)&=&\frac{t-t_m}{m(t_m-t_{m - 1})}+
1 ,\qquad t\in[t_{m},1] .
\end{eqnarray*}
In this case, we define $T_{4,\sigma}\dvtx[0,1]\times\RR\to[0,1]\times
\RR$
by $T_{4,\sigma}=(T_{1,\sigma}(t),x)$. Then,
$T_{1,\sigma}\dvtx[0,t_m]\to[0,1]$ and
$T_{4,\sigma}\dvtx[0,t_m]\times\RR\to[0,1]\times\RR$ are bijections
and we
let $T_{2,\sigma}\dvtx[0,1]\to[0,t_m]$ and $T_{3,\sigma}\dvtx[0,1]\times
\RR\to[0,t_m]\times\RR$ be their corresponding inverses.

Let $n\in\NN$. Recall \eqref{defVnk} and set
\begin{eqnarray*}
M_{1,n}&=&\sum_{k=1}^n
\varepsilon_{V_{n,k}} 1_{G_{1,n}(k/n)-G_{1,n}((k - 1)/n)\neq 0} ,\\
M_{2,n}&=&\sum_{0\le t\le1}
\varepsilon_{([tn]+1)/n} 1_{G_{2,n}(t)-G_{2,n}(t-)\neq0} .
\end{eqnarray*}
For $n\in\NN$ and $i=1,2$, it follows from the transformation theorem
that
%
\begin{eqnarray}\label{timechange1}
G_{i,n}\circ T_{i,M_{i,n}}(t)&=&\int_{[0,t]\times\RR}
(h_{n,i}\circ T_{M_{i,n}})^{1/2}(s-) z N_{i,n}^{T_{i+2,M_{i,n}}}(\mathrm{d}s,\mathrm{d}z),\nonumber \\
 h_{i,n}\circ
T_{i,M_{i,n}}(t)&=&h_0 + \int_{[0,t]} \beta - \alpha
(h_{i,n}\circ T_{i,M_{i,n}})(s-)\,\mathrm{d}T_{i,M_{i,n}}(s)\\
&&{} + \lambda\int_{[0,t]\times\RR}
(h_{i,n}\circ T_{i,M_{i,n}})(s-) z^2 N_n^{ T_{i+2,M_{i,n}}}(\mathrm{d}s,\mathrm{d}z)\nonumber
\end{eqnarray}
for all $t\in[0,1]$ and
$\theta=(h_0,\beta,\alpha,\gamma)\in[0,\infty)^4$.

Let $\theta=(h_0,\beta,\alpha,\gamma)\in[0,\infty)^4$. If
$h_0 = \beta = 0$, then it follows from \eqref{G1p}, \eqref{EEE3}
and \eqref{timechange1} that $h_{i,n}=h_{i,n}\circ T_{i,M_{i,n}}\equiv
0$, $i=1,2$, a.s. and thus
\[
\LLL_\theta(G_{i,n}) = \LLL_{\theta}(G_{i,n}\circ
T_{i,M_{i,n}}) = \varepsilon_0,\qquad  n\in\NN , i=1,2 .
\]
Otherwise, if $h_0+\beta>0$, it follows from \eqref{G1p} and
\eqref{EEE3} that $h_{i,n}(t)>0$ for all $t\in(0,1]$ a.s., $i=1,2$. In
this case, we have $M_{1,n} = N_{2,n}$, $M_{2,n} = N_{1,n}$,
$N_{1,n}^{T_{3,M_{1,n}}} = N_{n,2}$ and
$N_{2,n}^{T_{4,M_{2,n}}} = N_{n,1}$ and thus we get from
\eqref{timechange1} that both
\[
\LLL_\theta(G_{1,n}) = \LLL_{\theta}(G_{2,n}\circ
T_{2,M_{2,n}}) \quad\mbox{and}\quad
\LLL_\theta(G_{2,n})=\LLL_{\theta}(G_{1,n}\circ T_{1,M_{1,n}})
\]
for
$n\in\NN$. In other words, for all $n\in\NN$, there are Markov kernels
$K_{1,2,n}\dvtx D_1\times\DDD_1\to[0,1]$ and $K_{2,1,n}\dvtx D_1\times\DDD
_1\to
[0,1]$, not depending on $\theta\in[0,\infty)^4$, such that
$K_{1,2,n}\LLL_\theta(G_{2,n})=\LLL_\theta(G_{1,n})$ and
$K_{2,1,n}\LLL_\theta(G_{1,n}) = \LLL_\theta(G_{2,n})$ for all
$\theta\in[0,\infty)^4$. Hence,~$\EEE_{1,n}^{\star}$ is equivalent to~$\EEE_{2,n}^{\star}$ in deficiency by \eqref{defMarkov} for all
$n\in\NN$. This completes the proof of the theorem.

%
\section{\texorpdfstring{Proof of Theorem \protect\ref{th2a}}{Proof of Theorem 2.2}}
\vspace*{2pt}\label{proofB}
The proof of Theorem \ref{th2a} is split into two parts, reported in
Sections \ref{proofth2a} and \ref{SUBSECproofThm2aPartII},
respectively. We will need two additional results, which appear as
Lemmas \ref{Lemexpectbound} and \ref{condTVlem}, together with their respective proofs in
Sections \ref{SUBSEClemma51} and \ref{SUBSEClemma52}.
%
\vspace*{2pt}\subsection{\texorpdfstring{Proof of Theorem \protect\ref{th2a} (part I)}
{Proof of Theorem 2.2 (part I)}}\vspace*{2pt}\label{proofth2a}
Recall that Le Cam's distance is a pseudo-metric. In view of
\eqref{deflecamsmallerTV} and Theorem \ref{th2}, it thus suffices to
show \eqref{th2astronger}. For $n\in\NN$, let $Z_n=(Z_{n,k})_{1\le
k\le
n}$ be a random vector with a~distribution as in \eqref{Znthinn}.

First, we assume that
%
\begin{equation}\label{QnequalsQ}
Q_n = Q ,\qquad  n\in\NN .
\end{equation}
At the end of the proof we will relax this condition to $\|Q_n-Q\|\to
0$ as $n\to\infty$.

Let $N_n$ be as in \eqref{defNn} and set $\|N_n\|=N_n([0,1]\times
\RR)$, $n\in\NN$. Let $\Theta$ be as in the assertion of the theorem.
Suppose that
$H_{1,n}=H_n=(h_{0,1,n},\beta_{1,n},\alpha_{1,n},\lambda
_{1,n})\dvtx\Theta\to[0,\infty)^4
$ satisfies the assumptions of the theorem. Further, let
$H_{2,n}=(h_{0,2,n},\beta_{2,n},\alpha_{2,n},\lambda
_{2,n})=H_{n}^{(0)}\dvtx\Theta\to[0,\infty)^4
$ be defined by the identities in
\eqref{specialHn1} and \eqref{specialHn2}.

For $\theta\in\Theta$ and $i=1,2$, we define
$X_{i,n}=(X_{i,n}(k))_{1\le k\le n}$ by
%
\begin{eqnarray}\label{Xip}
X_{i,n}(k)&=&h^{1/2}_{i,n}(k - 1) Z_{n,k},\qquad
X_{i,n}(0)=0 ,\nonumber\\
h_{i,n}(k)&=&\beta_{i,n}(\theta)+
h_{i,n}(k - 1)  [\alpha_{i,n}(\theta)+\lambda_{i,n}(\theta)
 Z_{n,k}^2 ] ,\\
h_{i,n}(0)&=&h_{0,i,n}(\theta) ,\qquad  n\in\NN , 1 \le
k \le  n .\nonumber
\end{eqnarray}
Hence, $X_{1,n}$ corresponds to the GARCH processes $G_n$ as in the
theorem, and $X_{2,n}$ to the GARCH processes $G_n^{(0)}$ defined
directly after \eqref{specialHn2}. Let
%
\begin{equation}\label{defMnk}
M_{n,k} =  \Biggl\{\sigma=(\sigma
_l)_{1\le l\le k}\in\NN^k \dvtx  \sum_{l=1}^k\sigma_l\le n \Biggr\},\qquad 1 \le
k \le  n , n\in\NN.
\end{equation}
By employing the conventions
$0^0=1$ and $\sum_{l=k}^m=0$ for $m < k$, we set
%
\begin{eqnarray}\label{defhateta}
\eta_{i,n,1,l,\sigma}(\theta)&=&\beta_{i,n}(\theta)\sum
_{m=0}^{\sigma_{l + 1} - 1}[\alpha_{i,n}(\theta)]^m ,\nonumber\\
\eta_{i,n,2,l,\sigma}(\theta)&=&[\alpha_{i,n}(\theta)]^{\sigma
_{l + 1}} ,\\
\eta_{i,n,3,l,\sigma}(\theta)&=&\lambda_{i,n}(\theta) [\alpha
_{i,n}(\theta)]^{\sigma_{l + 1} - 1}\nonumber
\end{eqnarray}
for $\sigma=(\sigma_l)_{1\le l\le k}\in M_{n,k}$, $1 \le  k \le
n$, $0 \le  l \le  k - 1$, $i=1,2$ and $n\in\NN$.

Also, we recursively define functions from $\RR^k\to\RR$ by setting
%
\begin{eqnarray}\label{defhatg}
\hat g_{i,n,0,\sigma,\theta}&\equiv&h_{0,i,n}(\theta) \alpha
_{i,n}(\theta)^{\sigma_1-1}+\beta_{i,n}(\theta)
\sum_{m=0}^{\sigma_1-2}\alpha_{i,n}^m(\theta) ,\nonumber\\ [-8pt]\\ [-8pt]
\hat g_{i,n,l,\sigma,\theta}(y)&=&\eta_{i,n,1,l,\sigma}(\theta)+\eta
_{i,n,2,l,\sigma}(\theta)
\hat g_{i,n,l - 1,\sigma,\theta}(y)+\eta_{i,n,3,l,\sigma}(\theta)
 y_l^2\nonumber
\end{eqnarray}
for $y\in\RR^k$, $\sigma=(\sigma_l)_{1\le l\le k}\in M_{n,k}$,
$1 \le  k \le  n$, $0 \le  l \le  k - 1$, $i=1,2$ and
$n\in\NN$.

Let $n\in\NN$ and $1 \le  k \le  n$. On $\{\|N_n\|=k\}$, we
consider the following stopping times:
\[
\tau_{0}=0 ,\qquad \tau_{m}=\min\bigl\{\nu\in\{\tau_{m-1} + 1,\dots
,n\}\dvtx Z_{n,\nu}\neq0\bigr\},\qquad 1 \le
m \le  k .
\]
Using these stopping times, let $\Delta\tau=((\Delta
\tau_m)_{1\le m\le k})\in M_{n,k}$ be the random vector defined
componentwise by $\Delta\tau_m=\tau_m-\tau_{m - 1}$ for $1 \le
m \le  k$.

Let $i=1,2$, $n\in\NN$, $1 \le  k \le  n$ and $\theta\in
\Theta$. On
$\{\|N_n\|=0\},$ set $Y_{i,n}=0$, and otherwise,
%
\begin{equation}\label{defYipfam}
Y_{i,n} = (Y_{i,n}(l))_{1\le l\le\|N_n\|} = (X_{i,n}(\tau
_l))_{1\le
l\le\|N_n\|} .
\end{equation}
In the notation of \eqref{defhateta} and \eqref{defhatg}, $Y_{i,n}$
satisfies the following recursion on $\{\|N_n\|=k\}$:
%
\begin{eqnarray}\label{Yip}
Y_{i,n}(l) & = &
g^{1/2}_{i,n}(l - 1) Z_{n,\tau_l} ,\qquad
Y_{i,n}(0) = 0 ,\qquad 1 \le  l \le k ,\nonumber \\
g_{i,n}(l) & = & \eta_{i,n,1,l,\Delta\tau}(\theta)+
g_{i,n}(l - 1) \eta_{i,n,2,l,\Delta\tau}(\theta)\nonumber\\ [-8pt]\\ [-8pt]
&&{} +\eta_{i,n,3,l,\Delta\tau}(\theta)
 g_{i,n}(l - 1) Z_{n,\tau_l}^2 ,\qquad  1 \le  l \le
k - 1 ,\nonumber \\
g_{i,n}(0)&=&\hat g_{i,n,0,\Delta\tau,\theta} .\nonumber
\end{eqnarray}
Recall \eqref{defMnk}. For all $n\in\NN$, $1 \le  k \le  n$ and
$\sigma=(\sigma_l)_{1\le l\le k}\in M_{n,k}$, let
%
\begin{equation}\label{defANK}
A_{n,k,\sigma} =  \{
\|N_n\| = k , \Delta\tau = \sigma \} .
\end{equation}
For
future purposes, we collect some useful inequalities in the next lemma.

%
\subsection{\texorpdfstring{Statement and proof of Lemma \protect\ref{Lemexpectbound}}
{Statement and proof of Lemma 5.1}}\label{SUBSEClemma51}
\begin{lemma}\label{Lemexpectbound}
Suppose that $(\Theta,(H_n)_{n\in
\NN})$ satisfies the assumptions of Theorem \ref{th2a}.
Let $S\in(0,\infty)$ and suppose that $Q([-S,S])=1$.

There then exist some $C=C(S,\Theta)\in(1,\infty)$ and
$n_0=n_0(S,\Theta)\in\NN$ such that the following three inequalities
hold:
%
\begin{eqnarray} \label{hatg0diff}
|\hat
g_{1,n,0,\sigma,\theta}-\hat
g_{2,n,0,\sigma,\theta}|&\le& \frac Cn,\\\label{awayfromzero}
\hat g_{i,n,l,\sigma,\theta}(y)&\ge&C^{-1},\\\label{hatcondEgdiff}
E_\theta [ |\hat g_{1,n,l,\sigma,\theta}(Y_{1,n})- \hat
g_{2,n,l,\sigma,\theta}(Y_{1,n})|  | A_{n,k,\sigma} ] &\le&
\frac{C^k}n
\end{eqnarray}
for all $n \ge  n_0$, $1 \le  k \le  n$, $0 \le  l\le
k - 1$, $\sigma\in M_{n,k}$, $i=1,2$, $\theta\in\Theta$, $y\in
\RR^k$\vadjust{\goodbreak}
and $i=1,2$.
\end{lemma}
\begin{pf}
Let $(\Theta,(H_n)_{n\in\NN})$ be as in Theorem \ref{th2a}. First, note
that $(\Theta,(H_{1,n})_{n\in\NN})=(\Theta,(H_n)_{n\in\NN})$ satisfies
the assumptions of Theorem \ref{th2a}. Also, recall that
$(\Theta,(H_{2,n})_{n\in\NN})=(\Theta,(H^{(0)}_{n})_{n\in\NN})$ is
defined in \eqref{specialHn1} and \eqref{specialHn2}. In particular,
observe that $\alpha_{i,n}(\theta)\to1$ uniformly for all
$\theta=(h_0,\beta,\alpha,\lambda)\in\Theta$ as $n\to\infty$, $i=1,2$,
and thus there is an $n_1=n_1(\Theta)\in\NN$ satisfying
%
\begin{equation}\label{deltahochn}
\frac{\mathrm{e}^{-1}}2 \le[\alpha
_{i,n}(\theta)]^n=\exp\bigl(n\log\bigl[n+n\bigl(\alpha_{i,n}(\theta) -
1\bigr)\bigr]-n\log
n\bigr) \le 2 \mathrm{e}
\end{equation}
for all $n\ge n_1$, $i=1,2$ and
$\theta=(h_0,\beta,\alpha,\lambda)\in\Theta$.

It follows from our assumptions on $(\Theta,(H_n)_{n\in\NN})$ that
there exist $n_0=n_0(\Theta)\ge n_1$ and $C_1=C_1(\Theta)\in
(1,\infty)$
such that\vspace*{-2pt}
%
\begin{eqnarray}\label{lowerbound}
\hat g_{i,n,l,\sigma,\theta}(y)&\ge&
h_{0,i,n}(\theta)[\alpha_{i,n}(\theta)]^{-1+\sum_{m=1}^{l+1}\sigma_k}
\ge \frac{\mathrm{e}^{-1}}
2\frac{h_{0,i,n}(\theta)}{\alpha_{i,n}(\theta)}\nonumber\\ [-7pt]\\ [-7pt]
&\ge& \frac{\mathrm{e}^{-1}}
4\inf_{(h_0,\beta,\alpha,\lambda)\in\overline\Theta}h_{0}
 \ge C_1^{-1}\nonumber
\end{eqnarray}
and
%
\begin{eqnarray}\label{upperbound}
&&\max \biggl\{h_{0,i,n}(\theta) , \beta_{i,n}(\theta) ,
[\alpha_{i,n}(\theta)]^n ,
\frac{h_{0,i,n}(\theta)}{\alpha_{i,n}(\theta)}  \biggr\} \le
C_1 , \nonumber\\[1pt]
&&\max \{|h_{0,1,n}(\theta) - h_{0,2,n}(\theta)| ,
|\beta_{1,n}(\theta) - \beta_{2,n}(\theta)| ,\\[1pt]
&&\phantom{\max \{}|\alpha_{1,n}(\theta) - \alpha_{2,n}(\theta)| ,
|\lambda_{1,n}(\theta) - \lambda_{2,n}(\theta)| \} \le
\frac{C_1}n\nonumber
\end{eqnarray}
%
for all $n\ge n_0$, $1 \le  k \le  n$, $0 \le  l \le  k - 1$,
$\sigma=(\sigma_l)_{1\le l\le k}\in M_{n,k}$, $i=1,2$,
$\theta\in\Theta$ and\break $y\in\RR^k$.

Recall \eqref{defhateta} and \eqref{defhatg}. It follows from
\eqref{upperbound} that we have
%
\begin{eqnarray}\label{upperbound2}
\max \{\eta_{i,n,2,l,\sigma}(\theta) ,
\eta_{i,n,3,l,\sigma}(\theta) \}&\le& C_1^2 ,\nonumber\\[2pt]
\max \{\eta_{i,n,1,l,\sigma}(\theta) ,  \hat
g_{i,n,0,\sigma,\theta} \}&\le& (k + 1) C_1^2 ,\nonumber\\[2pt]
\max \{|[\alpha_{1,n}(\theta)]^m-[\alpha
_{1,n}(\theta)]^m| \}&\le&
\frac{C_1^2m} n ,\\[2pt]
\max \{|\eta_{1,n,j,l,\sigma}(\theta)-\eta
_{2,n,j,l,\sigma}(\theta)| \dvtx j=1,2,3 \}&\le&
\frac{(2\mathrm{e}^2 C_1^3)^k}n ,\nonumber\\[2pt]
|\hat g_{1,n,0,\sigma,\theta}-\hat
g_{2,n,0,\sigma,\theta}|&\le&\frac{(4 \mathrm{e}^2 C_1^3)^k}n\nonumber
\end{eqnarray}
for all $n\ge n_0$, $1 \le  k \le  n$, $0 \le  l\le  k - 1$,
$\sigma=(\sigma_l)_{1\le l\le k}\in M_{n,k}$, $i=1,2$, $m\in\NN_0$ and
$\theta\in\Theta$.

Recall \eqref{Yip}. Let $S>1$ be such that $Q([-S,S])=1$ and set
$C_2=C_2(S,\theta)=\mathrm{e}^2(1+S)^2C_1^4$ and $C_3=C_3(S,\theta)=S^2C_2$. It
follows from an induction and the inequalities in~\eqref{upperbound2}
that
%
\begin{eqnarray}\label{upperbound3}
E_\theta[g_{i,n}(l)|A_{n,k,\sigma}]
&\le&C_1^2(k + 1) + C_1^2(1 + S^2)E_\theta[g_{i,n}(l -
1)|A_{n,k,\sigma}]\nonumber \\ [-10pt]\\ [-10pt]
 &\le&(k + 1)\sum_{m=0}^l
(1 + S^2)^m C_1^{2(1+m)}\le C_2^k\nonumber
\end{eqnarray}
and thus that
%
\begin{equation}\label{upperbound4}
E_\theta[Y^2_{i,n}(l)|A_{n,k,\sigma}] \le
C_3^k
\end{equation}
for all $n\ge n_0$, $1 \le  k \le  n$, $0 \le  l \le  k - 1$,
$\sigma=(\sigma_l)_{1\le l\le k}\in M_{n,k}$, $i=1,2$ and
$\theta\in\Theta$.

Finally, let $C=C(S,\theta)=12 \mathrm{e}^3 C_1^3C_3$. By an induction, it
follows from \eqref{upperbound2}--\eqref{upperbound4}~that
\begin{eqnarray*}
&&E_\theta [  |\hat
g_{1,n,l,\sigma,\theta}(Y_{1,n})-
\hat g_{2,n,l,\sigma,\theta}(Y_{1,n}) | |A_{n,k,\sigma} ] \\[-2pt]
&&\quad\le C^k_3\sum_{j=1}^3|\eta_{1,n,j,l}(\theta)-\eta
_{2,n,j,l}(\theta)|
+ E_\theta [  |\hat g_{1,n,l - 1,\sigma,\theta}(Y_{1,n})-
\hat
g_{2,n,l - 1,\sigma,\theta}(Y_{1,n}) | |A_{n,k,\sigma} ]\\[-2pt]
&&\quad\le |\hat g_{1,n,0,\sigma,\theta}- \hat
g_{2,n,0,\sigma,\theta} |+C^k_4\sum_{l=1}^{k-1}\sum
_{j=1}^3|\eta_{1,n,j,l}(\theta)-\eta_{2,n,j,l}(\theta)| \le
\frac{C^k}n
\end{eqnarray*}
for all $n\ge n_0$, $1 \le  k \le  n$, $0 \le  l \le  k - 1$,
$\sigma=(\sigma_l)_{1\le l\le k}\in M_{n,k}$ and $\theta\in\Theta$.
This completes the proof in view of \eqref{lowerbound} and
\eqref{upperbound2}.
\end{pf}

%

\subsection{\texorpdfstring{Statement and proof of Lemma \protect\ref{condTVlem}}
{Statement and proof of Lemma 5.2}}\label{SUBSEClemma52}
We now provide an upper bound for conditional laws and their total
variation norms in the next lemma.
\begin{lemma}\label{condTVlem}
Suppose that $Q$ admits a Lebesgue
density $f$, where
$f$ is globally Lipschitz and has a compact support
$\overline{\{f>0\}}$.

If $(\Theta,(H_n)_{n\in\NN})$ satisfies the assumptions of
Theorem \ref{th2a}, then there exist $n_0=n_0(f,\Theta)\in\NN$ and
$C=C(f,\Theta)\in(0,\infty)$ such that
%
\begin{equation}\label{condTV}
\bigl\|\LLL_\theta (Y_{1,n}
|A_{n,k,\sigma} )-\LLL_\theta (Y_{2,n} |A_{n,k,\sigma} ) \bigr\| \le \frac
{C^k}n
\end{equation}
for all $\theta\in\Theta$, $n\ge n_0$, $1 \le
k \le  n$ and $\sigma\in M_{n,k}$.
\end{lemma}
\begin{pf}
By assumption, we have $f(x)=0$ for all $|x|\ge S$ and some $S>0$.
Hence, there are $n_0=n_0(f,\theta)\in\NN$ and
$C_1=C_1(f,\theta)\in(1,\infty)$ such that for $C$ replaced by $C_1$,
the assertion of Lemma \ref{Lemexpectbound} holds.

Let $n\ge n_0$, $i=1,2$, $\theta\in\Theta$, $1 \le  k \le  n$ and
$\sigma\in M_{n,k}$. Recall \eqref{defhatg}. In view of
\eqref{awayfromzero}, $\Psi_{i,n,\theta}\dvtx\RR^k\to\RR^k$ is
a\vadjust{\goodbreak}
well defined $C^\infty$ diffeomorphism, where
$\Psi_{i,n,\sigma,\theta}=(\psi_{i,n,l,\sigma,\theta})_{1\le l\le
k}\dvtx\break\RR^k\to\RR^k$ is defined by
%
\begin{eqnarray}\label{Phiip}
\psi_{i,n,l,\sigma,\theta}(y)=\frac{y_l}{\hat
g^{1/2}_{i,n,l - 1,\sigma,\theta}(y)}
\end{eqnarray}
for $y=(y_1,\dots,y_k)\in\RR^k$ and $1\le l\le k$. For all $n\ge n_0$,
$\theta\in\Theta$, $n\ge n_0$, $1 \le  k \le  n$ and $\sigma
\in
M_{n,k}$, we define
\[
\tilde f_{i,n,k,\sigma,\theta}(y) =
\prod_{l=1}^k \frac{f(\psi_{i,n,l,\sigma,\theta}(y))}{\hat
g^{1/2}_{i,n,l - 1,\sigma,\theta}(y)},\qquad  y\in\RR^k ,
i=1,2 .
\]
It follows from \eqref{defhatg}, \eqref{Yip} and \eqref{Phiip} that
$\tilde f_{i,n,k,\sigma,\theta}$ is a density of the probability
measure $\LLL_\theta (Y_{i,n} |A_{n,k,\sigma} )$ with respect
to the Lebesgue measure $\ell^{\otimes k}$ on $\BBB(\RR^k)$. In
particular, we must have
%
\begin{equation}\label{condequal}
 \bigl\|\LLL_\theta (Y_{1,n} |A_{n,k,\sigma} )
-\LLL_\theta (Y_{2,n} |A_{n,k,\sigma} ) \bigr\|  =
\frac
12\int_{\RR^k}|\tilde f_{1,n,k,\sigma,\theta}(y)-\tilde
f_{2,n,k,\sigma,\theta}(y)|\,\mathrm{d}y
\end{equation}
for all $\theta\in\Theta$, $n\ge n_0$, $1 \le  k \le  n$ and
$\sigma\in M_{n,k}$.

Suppose that $C_f\in(0,\infty)$ is a
global Lipschitz constant of $f$. By means of simple substitutions, for
all $\epsilon>0$ and $w, v\ge\epsilon$, we can observe that
\[
\frac
12 \int \biggl|\frac{f(x/v)}{v}-\frac{f(x/w)}{w} \biggr|\,\mathrm{d}x \le
\frac1
\epsilon (S^2 C_f+1) |v-w| .
\]
Consequently, for all $\epsilon>0$,
we find a $\kappa_1=\kappa_1(f,\epsilon)\in(1,\infty)$ such that
\[
\frac12
\int \biggl|\frac{f(x/v)}{v}-\frac{f(x/w)}{w} \biggr|\,\mathrm{d}x \le
\kappa_1(\epsilon)|v-w|,\qquad  v,w\ge\epsilon .
\]
In view of
\eqref{awayfromzero}, there thus exists some
$\kappa_2=\kappa_2(f,\Theta)\in(1,\infty)$ such that
%
\begin{eqnarray}\label{cinproof}
&&\frac12 \int \biggl|\frac{f(y_l/\hat
g^{1/2}_{1,n,l - 1,\sigma,\theta}(y))}{\hat
g_{1,n,l - 1,\sigma,\theta}^{1/2}(y)}- \frac{f(y_l/\hat
g^{1/2}_{2,n,l - 1,\sigma,\theta}(y))}{\hat
g^{1/2}_{2,n,l - 1,\sigma,\theta}(y)} \biggr|\,\mathrm{d}y_l\nonumber\\ [-8pt]\\ [-8pt]
&&\quad\le\kappa_2|\hat g_{1,n,{l - 1},\sigma,\theta}(y)-\hat
g_{2,n,l - 1,\sigma,\theta}(y)|\nonumber
\end{eqnarray}
for all $n\ge n_0$,
$1 \le  k \le  n$, $1 \le  l\le  k$, $\sigma\in M_{n,k}$,
$y\in\RR^{k}$ and $\theta\in\Theta$. By integrating over $y_k$, we get
from \eqref{cinproof} that
%
\begin{eqnarray}\label{cinproof2}
&&\frac12\int_{\RR^k}|\tilde f_{1,n,k,\sigma,\theta}(y)-\tilde
f_{2,n,k,\sigma,\theta}(y)|\,\mathrm{d}y\nonumber\\
&&\quad\le \kappa_2\int_{\RR^{k - 1}}\prod_{l=1}^{k-1}
\frac{f(\psi_{1,n,l,\sigma,\theta}(y))}{\hat
g^{1/2}_{1,n,l - 1,\sigma,\theta}(y)}  |\hat
g_{1,n,k - 1,\sigma,\theta}(y)-
\hat g_{2,n,k - 1,\sigma,\theta}(y) |\,\mathrm{d}y\\
&&\qquad{} +\frac12\int_{\RR^{k - 1}} \Biggl|\prod_{l=1}^{k-1}
\frac{f(\psi_{1,n,l,\sigma,\theta}(y))}{\hat
g^{1/2}_{1,n,l - 1,\sigma,\theta}(y)}- \prod_{l=1}^{k-1}
\frac{f(\psi_{2,n,l,\sigma,\theta}(y))}{\hat
g^{1/2}_{2,n,l - 1,\sigma,\theta}(y)} \Biggr|\,\mathrm{d}y\nonumber
\end{eqnarray}
for all $n\ge n_0$, $1 \le  k \le  n$, $\sigma\in M_{n,k}$ and
$\theta\in\Theta$. It follows from \eqref{hatcondEgdiff} that
%
\begin{eqnarray}\label{cinproof3}
&& \int_{\RR^{k - 1}}\prod_{l=1}^{k-1}
\frac{f(\psi_{1,n,l,\sigma,\theta}(y))}{g^{1/2}_{1,n,l - 1,\sigma
,\theta}(y)}
 |\hat g_{1,n,k - 1,\sigma,\theta}(y)- \hat g_{2,n,k -
1,\sigma,\theta}(y) |\,\mathrm{d}y\nonumber\\ [-8pt]\\ [-8pt]
&&\quad=E_\theta [  |\hat g_{1,n,k - 1,\theta}(Y_{1,n})-
\hat g_{2,n,k - 1,\theta}(Y_{1,n}) | |A_{n,k,\sigma}]\le
\frac{C_1^k}n\nonumber
\end{eqnarray}
for all $n\ge n_0$, $1 \le  k \le  n$, $\sigma\in M_{n,k}$ and
$\theta\in\Theta$.

Let $C=\mathrm{e} \kappa_2 C_1$. By induction, we thus get from
\eqref{hatg0diff} and \eqref{cinproof2}--\eqref{cinproof3} that
\[
\bigl \|\LLL_\theta (Y_{1,n} |A_{n,k,\sigma} )-\LLL
_\theta (Y_{2,n} |A_{n,k,\sigma} ) \bigr\|
\le\frac{C^k}{n} ,
\]
uniformly for all $n\,{\ge}\,n_0, 1\,{\le}\,k\,{\le}\,n, \sigma\,{\in}\,M_{n,k}$
and $\theta\,{\in}\,\Theta$. This completes the proof of the~lemma.\vspace*{-3pt}
\end{pf}

%

\subsection{\texorpdfstring{Proof of Theorem \protect\ref{th2a} (part II)}
{Proof of Theorem 2.2 (part II)}}\vspace*{-3pt}
\label{SUBSECproofThm2aPartII}
Let $f$ be a Lebesgue density of $Q$ and $\Theta$ be as in
Theorem \ref{th2a}. We denote the positive part of a function
$g\dvtx\RR\to\RR$ by $g^+$. Let $C^\infty_C$ be the space of infinitely
often continuously differentiable functions $g\dvtx\RR\to\RR$ with compact
support $\overline{\{g>0\}}$. As $C^\infty_C$ is dense in $L^1$, we
find a sequence of $g_m\in C_C^\infty$, $m\in\NN$, such that
$\int|g_m-f|\,\mathrm{d}\ell\to0$ as $m\to\infty$. It is immediate that both
$\int|g_m^+-f|\,\mathrm{d}\ell\to0$ and $\int g_m^+\,\mathrm{d}\ell\to1$ as $m\to
\infty$.
Without loss of generality, we may thus assume that $\int
g_m^+\,\mathrm{d}\ell>0$ for all $m\in\NN$. Then, $h_m:=g^+_m/\int g_m^+\,\mathrm{d}\ell$
defines a sequence of globally Lipschitz continuous probability
densities with a compact support $\overline{\{h_m>0\}}$ such that
$\int|h_m-f|\,\mathrm{d}\ell\to0$.

For $m\in\NN$, let $Z^{(m)}_n=(Z^{(m)}_{n,k})_{1\le k\le n}$ be a
random vector with distribution
\[
\LLL\bigl(Z^{(m)}_n\bigr)(B) =  \biggl((1 - p_n)\varepsilon_0(B)+p_n \int
_Bh_m\,\mathrm{d}\ell \biggr)^{\otimes n} ,
\]
with $B\in\BBB(\RR^n)$, $m,n\in\NN$, $1\le k\le n$. If we
replace $Z_{n,k}$ by $Z^{(m)}_{n,k}$ in \eqref{Xip}, then we get yet
another family of GARCH models, $X^{(m)}_{i,n}=(X^{(m)}_{i,n}(k))_{1\le
k\le n}$, say, indexed by $\theta\in\Theta$, $i=1,2$ and $m$,
$n\in\NN$.

It follows from the monotonicity theorem for Markov kernels and a
well-known upper bound for product measures (see \cite{Re93}, Lemma 1.4.2(i)
and page 23) that for all \mbox{$i=1,2$},
%
\begin{eqnarray}\label{truncatednoise}
\hspace*{-15pt}\sup_{\theta\in\Theta_0}\bigl\|\LLL_\theta(X_{i,n})-\LLL
_\theta\bigl(X^{(m)}_{i,n}\bigr)\bigr\|
&\le& \bigl\|\LLL(Z_{n})-\LLL\bigl(Z^{(m)}_n\bigr)\bigr\|\nonumber\\ [-8pt]\\ [-8pt]
\hspace*{-15pt}&\le&
n\bigl\|\LLL(Z_{n,1})-\LLL\bigl(Z^{(m)}_{n,1}\bigr)\bigr\|= \frac{n
p_n}2 \int|h_m-f|\,\mathrm{d}\ell .\nonumber
\end{eqnarray}\looseness=0
As $h_m$ is globally Lipschitz with a compact support
$\overline{\{h_m>0\}}$ for all $m\in\NN_0$, the assumptions of
Lemma \ref{condTVlem} hold. For all $m\in\NN$, there thus exist
$n_m\in\NN$ and $C_m=C(h_m,\Theta)\in(0,\infty)$ such that for all
$n\ge n_m$, we get, by conditioning and the monotonicity theorem for
Markov kernels, that
%
\begin{equation}\label{condTV1}
\sup_{\theta\in\Theta}\bigl\|\LLL_\theta
\bigl(X^{(m)}_{1,n}\bigr)-\LLL_\theta\bigl(X^{(m)}_{2,n}\bigr)\bigr\| \le \frac
1n  E \bigl[C_m^{\|N_n\|} \bigr]
\end{equation}
for $N_n$ as defined in \eqref{defNn}. Recall \eqref{Xip}. By combining
\eqref{truncatednoise} and \eqref{condTV1}, we get, from the triangle
inequality, that
\[
\sup_{\theta\in\Theta} \bigl\|\LLL_\theta(G_n)-\LLL_\theta
\bigl(G_n^{(0)}\bigr) \bigr\|
 \le  n p_n \int|h_m-f|\,\mathrm{d}\ell+\frac1n
 E \bigl[C_m^{\|N_n\|} \bigr]
\]
for all $m\in\NN$ and $n\ge n_m$. As
\eqref{lre} holds, we have $\lim_{n\to\infty}EC_m^{\|N_n\|}=
e^{\lambda(C_m-1)}$ and thus
\[
\limsup_{n\to\infty}\sup_{\theta\in\Theta} \bigl\|\LLL_\theta
(G_n)-\LLL_\theta\bigl(G_n^{(0)}\bigr) \bigr\|
\le\lambda \limsup_{m\to\infty}\int|h_m-f|\,\mathrm{d}\ell=0 ,
\]
giving
\eqref{th2astronger}. This completes the proof of Theorem \ref{th2a} in
the case where $Q_n=Q$ for all $n\in\NN$ (see \eqref{QnequalsQ}).

Now, assume that $Q_n\to Q$ in total variation norm as $n\to\infty$.
For $m\in\NN$, let $\widehat Z_n=(\widehat Z_{n,k})_{1\le k\le n}$ be a
random vector with distribution
\[
\LLL(\widehat Z_n) =  \bigl((1 - p_n)\varepsilon_0+p_n Q_n
\bigr)^{\otimes n},\qquad   n\in\NN .
\]
If we replace $Z_{n,k}$ by $\widehat Z_{n,k}$
in \eqref{Xip}, then we get the GARCH models in the assertion of the
theorem. We denote them by $\widehat X_{i,n}$, $n\in\NN$, $i=1,2$. By
the same argument as in~\eqref{truncatednoise}, we must have, for all
$i=1,2$ and $n\in\NN$ ,
\[
\sup_{\theta\in\Theta_0}\|\LLL_\theta(\widehat
{X}_{i,n})-\LLL_\theta(X_{i,n})\|  \le  np_n\|Q_n-Q\| .
\]
As the right-hand side tends to zero, this completes the proof of the
theorem.

%
\section{\texorpdfstring{Proof of Theorem \protect\ref{th1}}{Proof of Theorem 2.3}}\label{proofth1}
We first need to make some preparations. Let $Z=(Z_n)_{n\in\NN}$ and
$U=(U_n)_{n\in\NN}$ be independent sequences of i.i.d. random variables
such that $\LLL(Z_1)=Q$ with Lebesgue density~$f$ and $U_1$ is
uniformly distributed on $(0,1)$. For $d\in\NN$, we denote the order
statistics of $0,U_1,\ldots,U_d$ by $0=:U_{d,0}<U_{d,1}\le\cdots\le
U_{d,d}$. For each $n\in\NN$, let $\nu_n$ be a Poisson random variable
with parameter $\gamma_n>0$, independent of $Z$ and $U$.

In both \eqref{COG} and \eqref{hatEEE}, $N$ admits a representation $N=
\sum_{k=1}^{\nu_n} \varepsilon_{(U_{\nu_n,k},Z_k)}$ since $N$ is a~Poisson measure with intensity $\gamma_n\ell\otimes Q$. On
$\{\nu_n=0\},$ we let $\Delta U_{\nu_n}=\Delta G_{\nu_n}=\Delta
\widehat G_{\nu_n}=0$, whereas on $\{\nu_n>0\}$, we set
\begin{eqnarray*}
\Delta
U_{\nu_n}&=&(U_{\nu_n,k} - U_{\nu_n,k - 1})_{1\le
k\le\nu_n} ,\\[-2pt]
\Delta G_{\nu_n}&=&\bigl(G(U_{\nu_n,k})-G(U_{\nu_n,k}-)\bigr)_{1\le
k\le\nu_n} ,\\[-2pt]
\Delta\widehat G_{\nu_n}&=&\bigl(\widehat G(U_{\nu_n,k})-\widehat
G(U_{\nu_n,k}-)\bigr)_{1\le k\le\nu_n} .
\end{eqnarray*}
Let
$S_0=\RR^0=\{0\}$ and $\widetilde\RR = \bigcup_{d=0}^\infty\{d\}
\times S_d\times\RR^{d}$, where for $d\in\NN$, $S_d$ equals the set of
all $w=(w_1,\dots,w_d)'\in(0,1)^d$ such\vspace*{1pt} that $\sum_{i=1}^dw_i\le1$. We
endow $S_d$ and $\widetilde\RR$ with the Borel trace field $\BBB(S_d)$
($d\ge0$) and the $\sigma$-algebra $\widetilde\BBB$, respectively,
where $\widetilde\BBB$ is the set of all $B\subseteq\widetilde\RR$
such that $B\cap(\{d\}\times S_d\times\RR^{d})\in
\{\varnothing,\{d\}\}\otimes\BBB(S_d)\otimes\BBB(\RR^{d})$ for all
$d\in\NN_0$.

Since we have assumed that
$\Theta\subseteq(0,\infty)\times[0,\infty)^3$, and since $G$ and
$\widehat G$ always jump at the same time as $N$ does, all arrival
times are observed in full and thus $\EEE_{\gamma_n,Q}(\Theta)$ and~$\widehat\EEE_{\gamma_n,Q}(\Theta)$ are equivalent to $\FFF_n$ and
$\widehat\FFF_n$ in deficiency, respectively, in view of
\eqref{defMarkov}, where for all $n\in\NN$, we set
\begin{eqnarray*}
\FFF_n&=& (\widetilde\RR,\widetilde\BBB,(\LLL_{\theta}(\nu
_n,\Delta U_{\nu_n},\Delta G_{\nu_n}))_{\theta\in\Theta} ) ,\\[-2pt]
\widehat\FFF_n&=& (\widetilde\RR,\widetilde\BBB,(
\LLL_\theta(\nu_n,\Delta U_{\nu_n},\Delta\widehat
G_{\nu_n}))_{\theta\in\Theta} ).
\end{eqnarray*}
Let $\widehat w_0=0$ and for $d>0$, set $\widehat
w_d=(1/d,\dots,1/d)\in\RR^d$. Recall that
$\Theta\subseteq(0,\infty)\times[0,\infty)^3$ and pick $d\in\NN_0$,
$\theta=(h_0,\beta,\alpha,\lambda)\in\Theta$, $w=(w_1,\dots
,w_d)\in
S_d\cup\{\widehat w_d\}$. We define a diffeomorphism
$\Psi_{d,w,\theta}\dvtx\RR^d\to\RR^d$ as follows: if $d=0$, then let
$\Psi_{d,w,\theta}=0$; otherwise, if $d>0$,~let
\[
\Psi_{d,w,\theta}(z)= (h_{d,w,\theta,k}^{1/2}(z) z_k
)_{1\le k\le d},\qquad  z=(z_1,\dots,z_d)\in\RR^d ,
\]
where for $2\le k\le d$, we
recursively define
%
\begin{eqnarray}\label{defmap1}
h_{d,w,\theta,k}(z)&=&\frac{\beta}{\alpha}(1 - \mathrm{e}^{-\alpha w_k})+
\mathrm{e}^{-\alpha w_k}(1+\lambda z_{k-1}^2) h_{d,w,\theta,k-1}(z) ,\nonumber\\ [-9pt]\\ [-9pt]
h_{d,w,\theta,1}(z)&\equiv&h_{d,w,\theta,1}=\frac{\beta}{\alpha
}(1 - \mathrm{e}^{-\alpha
w_1})+ \mathrm{e}^{-\alpha w_1}h_{0},\nonumber
\end{eqnarray}
provided $\alpha>0$ and otherwise, if $\alpha=0$, we set
%
\begin{eqnarray}\label{defmap2}
h_{d,\theta,w,k}(z)&=&\beta w_k+ h_{d,\theta,w,k-1}(z)(1+\lambda
z_{k-1}^2) ,\nonumber \\ [-9pt]\\ [-9pt]
h_{d,\theta,w,1}(z)&\equiv&h_{d,w,\theta,1}=\beta w_1+
h_{0} .\nonumber
\end{eqnarray}
Let $f$ be a strictly positive Lebesgue density of $Q$ and set
\begin{eqnarray*}
&&\HHH_{d,\theta_1,\theta_2,w}(\zeta) \\[-2pt]
&&\quad= \int_{\RR^d}  ( |J_{\Psi^{-1}_{d,w,\theta
_1}}(x)| f^{\otimes
d}(\Psi^{-1}_{d,w,\theta_1}(x))  )^{ \zeta}  (
|J_{\Psi^{-1}_{d,w,\theta_2}}(x)| f^{\otimes
d}(\Psi^{-1}_{d,w,\theta_2}(x))  )^{ 1-\zeta}\,\mathrm{d}x
\end{eqnarray*}
for all $\theta_1,\theta_2\in\Theta$, $0 < \zeta < 1$, $w\in
S_d\cup\{\widehat{w}_d\}$.

To summarize, we have thus far shown that for all $n\in\NN$,
equivalence of $\EEE_{\gamma_n,Q}(\Theta)$ and
$\widehat\EEE_{\gamma_n,Q}(\Theta)$ in deficiency\vadjust{\goodbreak} is equivalent to
equivalence of $\FFF_{n}$ and $\widehat\FFF_{n}$ in deficiency. For the
remaining part, recall that the two experiments are equivalent in
deficiency if and only if their corresponding Hellinger transformations
are equal (see \cite{St85}, Corollary 53.8). By solving the
differential equations in \eqref{COG} and \eqref{hatEEE}, we thus
arrive at the following identity:\looseness=-1\vspace*{-1pt}
\[
\sum_{d=1}^\infty\frac{\gamma_n^d\mathrm{e}^{-\gamma_n}}{d!}
\HHH_{d,\theta_1,\theta_2,\widehat w_d}(\zeta) =  \sum
_{d=1}^\infty
\frac{\gamma_n^d\mathrm{e}^{-\gamma_n}}{d!}
\int_{S_d} \HHH_{d,\theta_1,\theta_2,w}(\zeta) \frac{\mathrm{d}w}{\ell
^{\otimes d}(S_d)}\vspace*{-2pt}
\]\looseness=0
for all $\theta_1,\theta_2\in\Theta$,
$0 < \zeta < 1$, $n\in\NN$.

In the last display, the functions are analytical in $\gamma_n$;
consequently, for all $d\in\NN$, $\theta_1,\theta_2\in\Theta$,
$0 < \zeta < 1$, we must have\vspace*{-1pt}
%
\begin{equation}\label{centralidentitynoneq}
\HHH_{d,\theta_1,\theta_2,\widehat w_d}(\zeta) =
\int_{S_d} \HHH_{d,\theta_1,\theta_2,w}(\zeta) \frac{\mathrm{d}w}{\ell
^{\otimes
d}(S_d)}.\vspace*{-2pt}
\end{equation}
Next, we return to the proof of the theorem. By our assumption, there
exists $\zeta_0\in(0,1)$ such that, with
$g_{f,\zeta_0}\dvtx(0,\infty)\to[0,1]$ as in \eqref{defgQz},
$g_{f, \zeta_0}$ is strictly increasing on $(0,1]$. As a result,
$h\mapsto g_{f, \zeta_0}(\sqrt{h})$ is strictly increasing on $(0,1]$.

For all $\theta_1,\theta_2\in\Theta$, define
$H_{\theta_1,\theta_2}\dvtx(0,1]\to(0,\infty)$ by
$H_{\theta_1,\theta_2}(w):=h_{1,w,\theta_2,1}(1)/h_{1,w,\theta_1,1}(1)$
for $0 < w \le  1$. In particular, taking $d=1$ and $\zeta=\zeta_0$
in \eqref{centralidentitynoneq}, we must have\vspace*{-1pt}
%
\begin{equation}\label{centralidentitynoneq2}
g_{f,\zeta_0} \bigl\{\sqrt{H_{\theta_1,\theta_2}(1)} \bigr\} =
\int_{(0,1)} g_{f,\zeta_0} \bigl\{\sqrt{H_{\theta_1,\theta
_2}(w)} \bigr\}\,\mathrm{d}w\vspace*{-2pt}
\end{equation}
for all $\theta_1,\theta_2\in\Theta$.

(i) and (ii)   For $i=1,2,$ let
$\theta_i=(h_{0,i},\beta_i,\alpha,\lambda)\in\Theta$. Then,\vspace*{-1pt}
\[
h^2_{1,w,\theta_1,1}(1) \mathrm{e}^{\alpha w} \frac{\mathrm{d}}{\mathrm{d}w}H_{\theta
_1,\theta_2}(w) = \beta_2 h_{0,1}-\beta_1 h_{0,2},\qquad 0< w\le
1 .\vspace*{-2pt}
\]
(Note that this formula extends to $\alpha=0$.) If
$\beta_1 = \beta_2 > 0$ and $h_{0,1} > h_{0,2}$, then
$H_{\theta_1,\theta_2}$ is strictly increasing with
$H_{\theta_1,\theta_2}(1) \le1$, contradicting
\eqref{centralidentitynoneq2} since $h\mapsto
g_{f, \zeta_0}(\sqrt{h})$ is strictly increasing on $(0,1]$. If
$h_{0,1} = h_{0,2}$ and $\beta_2 < \beta_1$, then
$H_{\theta_1,\theta_2}$ is strictly decreasing with
$H_{\theta_1,\theta_2}(0+)=1$, contradicting
\eqref{centralidentitynoneq2} since $h\mapsto
g_{f, \zeta_0}(\sqrt{h})$ is strictly increasing on $(0,1]$. Reversing
the roles of parameters by replacing $H_{\theta_1,\theta_2}$ with
$H_{\theta_2,\theta_1}$, the previous reasoning extends to the
remaining cases where either $\beta_1 = \beta_2 > 0$ and
$h_{0,1} < h_{0,2}$, or $h_{0,1} = h_{0,2}$ and
$\beta_2 > \beta_1$. This completes the proof
of (i) and (ii).

(iii)   If
$(h_{0},\beta,\alpha_1,\lambda),(h_{0},\beta,\alpha_2,\lambda)\in
\Theta$
and $\beta=0$, then we have
$H_{\theta_1,\theta_2}(w)=\mathrm{e}^{(\alpha_1-\alpha_2)w}$ for all
$w\in(0,1]$. (Note that this formula extends to $\alpha_1=0$ or
$\alpha_2=0$.) By the same arguments as in parts (i) and (ii), we get
from \eqref{centralidentitynoneq2} that
$\alpha_1=\alpha_2$.

(iv)   In view of (iii), we may assume that $\beta>0$. Contradicting
the hypothesis, we assume that $\alpha_2>0$. It follows from the strict
inequality $\mathrm{e}^x-1>x$, $x>0$, that\vspace*{-1pt}
\begin{eqnarray*}
(h_0+\beta w)^2 \frac{\mathrm{d}}{\mathrm{d}w}H_{\theta_1,\theta_2}(w)
&=&\mathrm{e}^{-\alpha_2w} \biggl\{w(\beta^2 - \alpha_2\beta
h_0)-h_0^2\alpha_2 - \frac{\beta^2}{\alpha_2}
(\mathrm{e}^{\alpha_2w} - 1) \biggr\}\\[-1pt]
&<&-\alpha_2 h_0 \mathrm{e}^{-\alpha_2 w}(h_0+w\beta) < 0\vspace*{-2pt}
\end{eqnarray*}
for all $w\in(0,1]$. Thus, $w\mapsto H_{\theta_1,\theta_2}(w)$ is
strictly decreasing on $(0,1]$ with $H_{\theta_1,\theta_2}(0+)=1$,
contradicting \eqref{centralidentitynoneq2}. Thus, we must have
$\alpha_2=0$.\vadjust{\goodbreak}

(v)   Let
$(h_{0},\beta,\alpha_1,\lambda),(h_{0},\beta,\alpha_2,\lambda)\in
\Theta$
with $\alpha_{2}>\alpha_{1}$. Without loss of generality, we may assume
that $\beta > 0$. First, assume that $\beta/\alpha_2\le h_0\le
\beta/\alpha_1$. Then, $\beta-\alpha_1 h_0\ge0$ and $\beta-\alpha_2
h_0\le0$. Note that we cannot simultaneously have that
$\beta - \alpha_1 h_0=\beta - \alpha_2 h_0= 0$ such that
\[
h_{1,w,\theta_1,1}^2 \frac{\mathrm{d}}{\mathrm{d}w}H_{\theta_1,\theta_2}(w) =
(\beta - \alpha_2h_0)\mathrm{e}^{-\alpha_2 w}h_{1,w,\theta_1,1}
-(\beta - \alpha_1h_0)\mathrm{e}^{-\alpha_1 w}h_{1,w,\theta_2,1}<0
\]
for all $0<w\le1$. Consequently, $H_{\theta_1,\theta_2}$ is strictly
decreasing with $H_{\theta_1,\theta_2}(0+)=1$,
contradicting \eqref{centralidentitynoneq2}. Second, let
$h_0<\beta/\alpha_2$ and set
\[
\psi(w):=
(\beta-\alpha_2h_0)h_{1,w,\theta_1,1}
-(\beta-\alpha_1h_0)\mathrm{e}^{-(\alpha_1-\alpha_2) w}h_{1,w,\theta_2,1},\qquad
 0<w\le1 .
\]
As we have $\alpha_2>\alpha_1$ and $h_0<\beta/\alpha_2$, we must have
that $\beta - \alpha_1h_0>\beta - \alpha_2h_0>0$ such~that
\[
\psi'(w)=(\alpha_1 - \alpha_2)(\beta - \alpha_1h_0)\mathrm{e}^{-(\alpha
_1 - \alpha_2) w}h_{d,w,\theta_2,1}<0,\qquad
 0<w\le1 .
\]
Note that $\psi(0+)=(\alpha_1 - \alpha_2)
h_0^2<0$ and thus $\psi(w)<0$ for all $0<w\le1$. Since\break $\mathrm{e}^{\alpha_2w
}h_{d,w,\theta_1,1}^2 \frac{\mathrm{d}}{\mathrm{d}w}H_{\theta_1,\theta_2}(w)=\psi(w)<0$
for all $0 < w \le  1$, $H_{\theta_1,\theta_2}$ is strictly
decreasing with $H_{\theta_1,\theta_2}(0+)=1$,
contradicting \eqref{centralidentitynoneq2}. This completes the proof
of (v).\

\vspace*{-3pt}\section{\texorpdfstring{Proofs of the results in Section \protect\ref{secComplete}}
{Proofs of the results in Section 2.4}}\vspace*{-3pt}\label{SUBSEC4.9}
This section contains the proofs of Propositions \ref{prop1}--\ref{prop3}.
%
\vspace*{-3pt}\subsection{\texorpdfstring{Proof of Proposition \protect\ref{prop1}}
{Proof of Proposition 2.1}}\vspace*{-3pt}\label{proofprop1} For
$f\in D_d[0,1]$, we write $\Delta f=f(t)-f(t-)$,
$0\le t\le1$, with the convention that $\Delta f(0)=0$. With the usual
convention $\inf\varnothing=\infty$, we define
$T(f)=\inf\{t\in[0,1]\dvtx\Delta f_1(t)\neq0\}\wedge1$ for all
$f=(f_1,f_2)\in D_2$. Let $S$ be the set of all functions $f\in D_2$
with $T(f)\in(0,1)$. Let $D''_{0}\subseteq D_2[0,1]$ be the set of all
functions $f=(f_1,f_2)$ such that the right-hand derivatives~$f'_1(0+)$
and $f''_2(0+)$ exist in $\RR$. Further, let $D''_{0,T}\subseteq S\cap
D_0''$ be the set of all functions $f=(f_1,f_2)$ such that the
right-hand derivatives $f'_2(T(f)+)$ and $f''_2(T(f)+)$ exist in $\RR$.

Let $f\in D_2$ with $T=T(f)$. If $f\in D''_{0,T}$ and $f'_2(0+)\neq0$,
then we set
\[
X(f) =  \biggl(|f_2(0)|, \biggl|\frac{(f_2'(0+))^2 -
f_2(0)f_2''(0+)}{f'_2(0+)} \biggr|,
 \biggl|\frac{f_2''(0+)}{f'_2(0+)} \biggr|,\frac{|\Delta
f_2(T)|}{(\Delta
f_1(T))^2} \biggr) .
\]
If $f\in D''_{0,T}$, $f'_2(0+)= 0$ and
$f'_2(T+)\neq0$, then we set
\[
X(f) =  \biggl(|f_2(0)|, \biggl|\frac{(f_2'(T+))^2 -
f_2(T)f_2''(T+)}{f'_2(T+)} \biggr|,
 \biggl|\frac{f_2''(T+)}{f'_2(T+)} \biggr|,\frac{|\Delta
f_2(T)|}{(\Delta
f_1(T))^2} \biggr) .
\]
If $f\in D''_{0,T}$, $f'(0+)=f'_2(T+)=0$ and
$\Delta f_2(T)\neq0$, then we set
\[
X(f) =  \biggl(|f_2(0)|,0,0
,\frac{|\Delta f_2(T)|}{(\Delta f_1(T))^2} \biggr) .\vadjust{\goodbreak}
\]
If $f\in
D''_{0,T}$, $f'(0+)=f'_2(T+)=0$ and $\Delta f_2(T)= 0$, then we set
$X(f)=(|f_2(0)|,\infty,\infty,0)$. If $f\in D''_{0}\backslash S$ and
$f'_2(0+)\neq0$, then we define
\[
X(f) =  \biggl(|f_2(0)|,  \biggl|\frac{(f_2'(0+))^2 -
f_2(0)f_2''(0+)}{f'_2(0+)} \biggr|,
 \biggl|\frac{f_2''(0+)}{f'_2(0+)} \biggr|,\infty \biggr) .
\]
For the
remaining cases, we set
$X(f)=  (|f_2(0)|,\infty,\infty,\infty )$. Then,
$X\dvtx D_2\to[0,\infty]^4$ is a $\DDD_2$-$\BBB([0,\infty]^4)$-measurable
mapping. Since $Q(\{0\})=0$, it follows from \eqref{COG} that
$\LLL_\theta^X((G,h))=Q_{\theta}$ for all $\theta\in[0,\infty)^4$ and
thus $\delta(\EEE_{h},\FFF)=0$ by \eqref{defMarkov}, where $\FFF$ is
the experiment as defined in the assertion of the proposition.

Next, we show that $\delta(\FFF,\EEE_h)=0$. To this end, we define
$\xi=(\xi_1,\dots,\xi_3)\dvtx[0,\infty]^4\to[0,\infty)^3$ as
follows. Let
$\omega=(\omega_1,\dots,\omega_4)\in[0,\infty]^4$. If
$(\omega_1,\dots,\omega_3)\in[0,\infty)^3$, then we set
$\xi(\omega)=(\omega_1,\omega_2,\omega_3)$; if $\omega_1\in
[0,\infty)$
and either $\omega_2=\infty$ or $\omega_3=\infty$, then we set
$\xi(\omega)=(\omega_1,0,0)$; otherwise, we set $\xi(\omega)=0$.

In the notation of the \hyperref[s1]{Introduction}, we define $\Psi:
[0,\infty)^3\times\MM_2\to D_2$, where, for $0\le t\le1$,
$\omega=(\omega_1,\omega_2,\omega_3)\in[0,\infty)^3$ and
$\sigma\in\MM_2$, $(f_1(t),f_2(t))=\Psi[\omega,\sigma](t)$ is defined
to be the unique solution of the system of the following integral
equations:
%
\begin{eqnarray}\label{COGth4}
\hspace*{-15pt}f_1(t) & =  & \int_{[0,t]\times\RR^2}
f_2^{1/2}(s-)z_1 \sigma(\mathrm{d}s,\mathrm{d}z_1,\mathrm{d}z_2) ,\nonumber \\ [-8pt]\\ [-8pt]
\hspace*{-15pt}f_2(t) &  =  & \omega_1 + \int_{[0,t]}\bigl(\omega_2-\omega_3
f_2(s-)\bigr)\,\mathrm{d}s + \int_{[0,t]\times\RR\times(0,\infty)}
f_2(s-)z_2^2 \sigma(\mathrm{d}s,\mathrm{d}z_1,\mathrm{d}z_2) .\nonumber
\end{eqnarray}
Clearly, $\Psi$ is
$(\BBB([0,\infty)^3)\otimes\MMM_2)/\DDD_2$-measurable and thus defines
a deterministic Markov kernel
$K_2\dvtx([0,\infty)^3\times\MM_2)\times\DDD_2\to[0,1]$.

Let $\nu_0$ be the zero measure on $\BBB([0,1]\times\RR^2)$. For
$\lambda\ge0$, let $M_{\lambda}$ be a Poisson measure on
$[0,1]\times
\RR^2$ with the intensity measure $\gamma\ell\otimes
\LLL(Z,\lambda^{1/2} Z)$, where $\LLL(Z)=Q$ and $\gamma>0$ is the
intensity parameter of $N$ in \eqref{COG}. Consider the Markov kernel
$K_1\dvtx[0,\infty]^4\times(\BBB([0,\infty)^3)\otimes\MMM_2)\to[0,1]$
defined by
\[
K_1 [(\omega_1,\omega_2,\omega_3,\omega_4), \cdot  ] =
\varepsilon_{\xi(\omega)} \otimes  \cases{
\varepsilon_{\nu_0} ,&\quad$\omega_4=\infty$,\cr
\LLL(M_{\omega_4}|M_{\omega_4}\neq
\nu_0) ,&\quad$\omega_4 < \infty$.
}
\]
Observe that $K_2K_1Q_{\theta}=\LLL_{\theta}(G,h)$ for all
$\theta\in[0,\infty)^4$, in view of \eqref{COG}. Hence,
$\delta(\FFF,\EEE_h)=0$, by~\eqref{defMarkov}.

To summarize, we have shown that $\EEE_h$ is equivalent to $\FFF$ in
deficiency. By means of similar arguments, we can show that
$\Delta(\FFF,\widehat\EEE_h)=0$.
%
\subsection{\texorpdfstring{Proof of Proposition \protect\ref{prop2}}{Proof of Proposition 2.2}}\label{proofprop2}
(i) Let
$H_n=H_n^{(0)}\dvtx[0,\infty)^4\to[0,\infty)^4$ be as defined in
\eqref{specialHn1}--\eqref{specialHn2} and define $\bar
H_n\dvtx[0,\infty)^3\to
M:=\{(x_1,x_2,x_3)\in[0,\infty)^2\times(0,1] \dvtx x_1\ge x_2\}$ by
\[
\bar
H_n(h_0,\beta,\alpha)= (h_{0,n}(h_0,\beta,\alpha,0),\beta
_n(h_0,\beta,\alpha,0),\alpha_n(h_0,\beta,\alpha,0) )
 ,
\]
$h_0,\beta,\alpha\in[0,\infty)$. Then, $H_n\dvtx[0,\infty)^3\to
M\times[0,\infty)$ and $\bar H_n\dvtx[0,\infty)^3\to M$ are both bijections
with inverse functions $H_n^{-1}\dvtx M\times[0,\infty)\to[0,\infty)^4$
and $\bar H_n^{-1}\dvtx M\to[0,\infty)^3$, respectively. Define
$\widetilde
H_n\dvtx \RR^3\to[0,\infty)^3$ and $\widehat H_n\dvtx\RR^4\to[0,\infty)^4$ by
$\widetilde H_n(x_1,x_2,x_3)=\bar H_n^{-1}(|x_1|\vee
|x_2|,|x_2|,|x_3|\wedge1)$ and $\widehat
H_n(x_1,x_2, x_3,x_4)=H_n^{-1}(|x_1|\vee|x_2|,|x_2|,|x_3|\wedge
1,|x_4|)$ for $x_1,x_2,\break x_3,x_4\in\RR$ with $x_3\neq0$.

In the sequel, we write $x=(x(k))_{0\le k\le n}$ for a generic element
of $\RR^{n + 1}$. Fix $n\ge5$. Let $M_{0,n}\subseteq
[\RR^{n + 1}]^2$ be the set of all $(x,y)$ such that both $y(0)\neq
y(1)$ and $y(1)\neq y(2)$. By employing the convention
$\inf\varnothing=\infty$, we define
$T_n\dvtx[\RR^{n + 1}]^2\to\{1,\dots,n + 1\}$ by
\[
T_n(x,y) = \inf\{1 \le  k \le  n\dvtx x(k)\neq
x(k - 1)\}\wedge1,\qquad  (x,y)\in[\RR^{n + 1}]^2 .
\]
Let $S_n$
be the set of all $(x,y)\in[\RR^{n + 1}]^2$ with $3 \le
T(x,y) \le  n - 2$ such that $x(T)=x(T + 1)=x(T + 2)$. Consider
the subset $M_{T,n}\subseteq S_n$ of all $(x,y)\in[\RR^{n + 1}]^2$
such that both $y(T)\neq y(T + 1)$ and $y(T + 1)\neq y(T + 2)$
are satisfied.

For all $n\ge5$, we define a mapping
$X_n\dvtx[\RR^{n + 1}]^2\to[0,\infty]^4$ as follows: fix
$(x,y)\in[\RR^{n + 1}]^2$ and set $T=T_n(x,y)$. If $(x,y)\in
S_n\cap
M_{0,n}$, then set
\begin{eqnarray*}
X_n(x,y) &=&\widehat H_n \biggl(y(0),
\frac{y(1)^2 - y(0)y(2)}{y(1) - y(0)} ,
\frac{y(2) - y(1)}{y(1) - y(0)},\\[3pt]
&&\hphantom{\widehat H_n \biggl(} \frac{y(T)}{[x(T) - x(T  -  1)]^2} - \frac{y(1)^2-y(0)y(2)+
y(T - 1)[y(2) - y(1)] }{[y(1) - y(0)][x(T) - x(T  -  1)]^2}
 \biggr) .
\end{eqnarray*}
If $(x,y)\in M_{T,n}\backslash M_{0,n}$, then
set
\begin{eqnarray*}
&&X_n(x,y)\\[3pt]
&&\quad =\widehat H_n \biggl(y(0) ,
\frac{y(T + 1)^2 - y(T)y(T + 2)}{y(T + 1) - y(T)} ,
\frac{y(T + 2) - y(T + 1)}{y(T + 1) - y(T)},\frac{ y(T)}{[x(T) - x(T  -  1)]^2}\\[3pt]
&&\phantom{\quad =\widehat H_n \biggl(}{}-\frac{y(T + 1)^2-y(T)y(T + 2)+
y(T - 1)[y(T + 2) - y(T + 1)]
}{[y(T + 1) - y(T)][x(T) - x(T  -  1)]^2}
 \biggr).
\end{eqnarray*}
If $(x,y)\in S_{n}\backslash(M_{0,n}\cup
M_{T,n})$ and $y(T)\neq y(T - 1)$, then set
\[
X_n(x,y)
= \biggl(y(0),0,0, \frac{|y(T) -
y(T - 1)|}{(x(T) - x(T  -  1))^2}  \biggr) .
\]
If
$(x,y)\in S_{n}\backslash(M_{0,n}\cup M_{T,n})$ and $y(T)=
y(T - 1)$, then set $X_n(x,y)=(|y(0)|,\infty,\infty,0)$. If
$(x,y)\in
M_{0,n}\backslash S_n$ and $T=n + 1$, then set
\[
X_n(x,y) = \biggl(\widetilde H_n \biggl[y(0),
\frac{y(1)^2 - y(0)y(2)}{y(1) - y(0)},
\frac{y(2) - y(1)}{y(1) - y(0)} \biggr] , \infty \biggr).
\]
Otherwise, set $X_n(x,y) = (|y(0)|,\infty,\infty,\infty)$.

Recall that both $G_n=(G_{n,k})_{0\le k\le n}$ and $h_n=(h_{n,k})_{0\le
k\le n}$ are defined by \eqref{G1} via~\eqref{specialHn1} and \eqref{specialHn2}. For $n\ge5$, the mapping
$X_n\dvtx[\RR^{n + 1}]^2\to[0,\infty]^4$ is well defined and
$\BBB([\RR^{n + 1}]^2)/ \BBB([0,\infty]^4)$-measurable. Recall that
$Q_n(\{0\})=0$ for all $n\in\NN$ and thus\looseness=1
\[
\LLL_\theta^{X_n}(G_n,h_n) =  \cases{
q_{1,n}\varepsilon_{(h_0,\beta,\alpha,\infty)} +
q_{2,n} \varepsilon_{\theta}\vspace*{2pt}\cr
\quad {}+(1  -  q_{1,n} - q_{2,n})\varepsilon_{(h_{0,n}(\theta
),\infty,\infty,\infty)} ,&\quad $\theta\notin\Theta_{e}$,\vspace*{2pt}\cr
(1  -  q_{2,n})
\varepsilon_{(h_{0},\infty,\infty,\infty)} + q_{2,n}
\varepsilon_{\theta} ,&\quad $\theta\in\Theta_{e} , h_0 > 0 ,
\lambda > 0$,\vspace*{2pt}\cr
(1  -  q_{2,n})
\varepsilon_{(h_{0},\infty,\infty,\infty)} + q_{2,n}
\varepsilon_{(h_{0},\infty,\infty,0)} ,&\quad $\theta\in\Theta_{e} ,
h_0 > 0 , \lambda = 0$,\vspace*{2pt}\cr
\varepsilon_{(0,\infty,\infty,\infty)} ,&\quad $\theta\in\Theta_{e}
,h_0 = 0$
}
\]\looseness=0
for all $n\ge5$, $\theta=(h_0,\beta,\alpha,\lambda)\in[0,\infty)^4$,
where we set $q_{1,n}=(1  -  p_n)^n$ and
$q_{2,n}=(1  -  p_n)^2[1 - p_n- p_n (1  -  p_n)][1 - (1
 -  p_n)^{n-4}]$.

On the other hand, define a mapping
$\xi_n=(\xi_{1,n},\dots,\xi_{4,n})\dvtx[0,\infty]^4\to[0,\infty)^4$ as
follows. Let $\omega=(\omega_1,\dots,\omega_4)\in[0,\infty]^4$. If
$\omega\in[0,\infty)^4$, then set $\xi_n(\omega)=H_n(\omega)$. If
$\omega\in[0,\infty)^3\times\{\infty\}$, then set $\xi_n(\omega
)=(\bar
H_n(\omega_1,\omega_2,\omega_3),0)$. If
$\omega\in[0,\infty)\times(\{\infty\}\times[0,\infty]\cup
[0,\infty]\times\{\infty\})\times[0,\infty)$, then set
$\xi_n(\omega)=(\omega_1,0,1,\omega_4)$. If
$\omega\in[0,\infty)\times(\{\infty\}\times[0,\infty]\cup
[0,\infty]\times\{\infty\})\times\{\infty\}$, then set
$\xi_n(\omega)=(\omega_1,0,1,0)$. Otherwise, set $\xi(\omega)=0$.
Define a Markov kernel
$K_{1,n}\dvtx[0,\infty]^4\times\BBB([0,\infty)^3\times
[\RR^{n}]^2)\to[0,1]$ by
\begin{eqnarray*}
K_{1,n} [\omega, \cdot  ]&=&
\varepsilon_{(\xi_{n,1}(\omega), \xi_{n,2}(\omega), \xi
_{n,3}(\omega))} \\[2pt]
&&\otimes  \cases{
\varepsilon_0 ,&\quad$\omega_4=\infty$,\vspace*{2pt}\cr
\LLL\bigl((Z_{n,k})_k,(\xi_4(\omega)Z^2_{n,k})_k|(Z_{n,k})_k\neq
0\bigr) ,&\quad $\omega_4 < \infty$
}
\end{eqnarray*}
for $\omega=(\omega_1,\omega_2,\omega_3,\omega_4)\in[0,\infty]^4$,
where $Z_n=(Z_{n,k})_k$ is the random vector with the distribution as
specified by \eqref{Znthinn}.

Also, let $K_{2,n}\dvtx[0,\infty)^3\times[\RR^{n}]^2\times\BBB(
[\RR^{n + 1}]^2)\to[0,1]$ be the Markov kernel defined by the
deterministic mapping $(\xi_1,\xi_2,\xi_3,z_1,z_2)\mapsto(x,y)$, where
we recursively set $x(0)=0$ and $y(0) = \xi_1$, and for $1\le k\le
n$,
\[
x(k)=x(k - 1)+y^{1/2}(k - 1) z_1(k),\qquad
y(k)=\xi_2+y(k - 1)\bigl(\xi_3+ z_2(k)\bigr) .
\]
For $n\ge5$, let
$\FFF_n=([0,\infty]^4,\BBB([0,\infty]^4),(\LLL_\theta^{
X_n}(G_n,h_n))_{\theta\in[0,\infty)^4})$. By construction, we have
$\delta(\EEE_{h,n},\FFF_n)=0$, by means of \eqref{defMarkov}. For all
$n\ge5$, observe that
\begin{eqnarray*}
\delta(\FFF_n,\EEE_{h,n}) & \le&
\sup_{\theta\in[0,\infty)^3}\|\LLL_\theta
(G_n,h_n)-K_{2,n}K_{1,n}\LLL^{X_n}_\theta(G_n,h_n)\| \\[2pt]
& \le& |1-q_{1,n}-q_{2,n}|+|1-(1-p_n)^n-q_{2,n}|
\end{eqnarray*}
%
and thus $\EEE_{h,n}$ is strongly asymptotically equivalent to $\FFF_n$
as $n\to\infty$, by means of \eqref{defMarkov} and \eqref{lre}. By
\eqref{deflecamsmallerTV}, $\FFF_n$ converges strongly to the
experiment $\FFF$ in the assertion of
Proposition \ref{prop1}, completing the proof of (i).

(ii)   This follows from the same arguments as in (i).
%
\subsection{\texorpdfstring{Proof of Proposition \protect\ref{prop3}}{Proof of Proposition 2.3}}\label{proofprop3}
(i)   Define $X,X_n\dvtx[0,\infty]^4\to[0,\infty]^4$ as follows. If
$\omega=(\omega_1,\dots,\omega_4)\in[0,\infty)^3\times\{\infty\}
$ such
that $\omega_1\omega_3=\omega_2$, then set
$X(\omega)=(\omega_1,\infty,\infty,\infty)$; otherwise, set
$X(\omega)=\omega$. If $\omega=(\omega_1,\dots,\omega_4)\in
[0,\infty)^3\times\{\infty\}$ such that $\omega_1
n(1-\mathrm{e}^{-\omega_3/n})=\omega_2$, then set
$X_n(\omega)=(\omega_1,\infty,\infty,\infty)$; otherwise, set
$X_n(\omega)=\omega$, $n\in\NN$.

By definition, the deficiency is non-decreasing in the parameter set
with respect to set inclusions. Further, we have $\widehat
Q^X_{\theta}=Q_{\theta}$ and $\widehat Q^{X_n}_{\theta}= Q_{\theta,n}$
for all $n\in\NN$ and thus, by~\eqref{defMarkov},
$\delta(\widehat\FFF(\Theta),\FFF(\Theta))\le
\delta(\widehat\FFF,\FFF)=0$ and
$\delta(\widehat\FFF(\Theta),\FFF_n(\Theta))\le
\delta(\widehat\FFF,\FFF_n)=0$
for all $n\in\NN$, completing the proof of (i).

(ii)   First, assume that $\Theta$ satisfies
\eqref{hsgleichbetaalphagleich} for all $x>0$. Without loss of
generality, we may assume that $\Theta\subseteq[0,\infty)^4$ is a
finite set (see \cite{St85}, Theorem 51.4). Define $\Omega_\Theta$ to
be the set of all $\omega=(\omega_1,\omega_2,\omega_3,\omega_4)
\in(0,\infty)\times\{\infty\}^2\times\{0,\infty\}$ such that $(
\omega_1,\beta,\alpha,\lambda)\in(\Theta\cap\Theta_e)\backslash
\widehat\Theta_e$
for some $(\beta,\alpha,\lambda)\in[0,\infty)^3$. If
$\omega=(\omega_1,\omega_2,\omega_3,\omega_4)\in\Omega_\Theta$,
then it
follows from \eqref{hsgleichbetaalphagleich} that the corresponding
pair $(\beta,\alpha)=(\beta(\omega_1),\alpha(\omega_1))\in
[0,\infty)^2$
is uniquely determined by $\omega_1$. Hence, we may define a mapping
$Y\dvtx[0,\infty]^4\to[0,\infty]^4$ as follows: if
$\omega=(\omega_1,\omega_2,\omega_3,\omega_4)\in\Omega_\Theta$,
then we
set $Y(\omega)=( \omega_1,\beta(\omega_1),\alpha(\omega_1),\omega_4)$;
otherwise, if $\omega\in[0,\infty]^4\backslash\Omega_\Theta$,
then we
set $Y(\omega)=\omega$. As $\Theta$, and thus $\Omega_\Theta$, is a
finite set, the mapping $Y$ is
$\BBB([0,\infty]^4)/\BBB([0,\infty]^4)$-measurable. In view of
\eqref{hsgleichbetaalphagleich}, note that $Q^Y_\theta=\widehat
Q_\theta$ for all $\theta\in\Theta$ and thus
$\delta(\FFF(\Theta),\widehat\FFF(\Theta))=0$, by~\eqref{defMarkov}.

Second, assume that \eqref{hsgleichbetaalphagleich} is violated. There
then exist some $h_0 > 0$,
$\theta_1=(h_{0},\beta_1,\alpha_2,\lambda_1)$
$\in\Theta\cap\Theta_{e}\cap\widehat\Theta_{e}^C$ and
$\theta_2=(h_{0},\beta_2,\alpha_2,\lambda_2)\in\Theta\cap\Theta_{e}$
such that $(\beta_1,\alpha_{1})\neq(\beta_2,\alpha_{2})$.

Consider $\Theta_{0}=\{\theta_{1},\theta_{2}\}$ and the decision space
$D=\{(\beta_1,\alpha_1),(\beta_2,\alpha_2)\}$, endowed with the
discrete topology. For $\theta=(h_0,\beta,\alpha,\lambda)\in\Theta$,
consider (continuous and bounded) loss functions $W_\theta\dvtx D\to\RR$,
where for $x=(x_1,\dots,x_4)\in[0,\infty]^4$, we set
$W_{\theta}(x)=1-1_{\{(\beta,\alpha)\}}(x_2,x_3)$. Further, we
define a
Markov kernel $\widehat\rho\dvtx[0,\infty]^4\times\BBB(D)\to[0,1]$, where
for $x\in[0,\infty]^4$ and $B\in\BBB(D)$, we set
\[
\widehat\rho(x,B) =  \cases{
\epsilon_{(\beta_1,\alpha_1)}(B) ,&\quad if  $x\in(0,\infty)\times
\{\beta_1\}\times\{\alpha_1\}\times[0,\infty)$,\cr
\epsilon_{(\beta_2,\alpha_2)}(B) ,&\quad otherwise.
}
\]
We then have $\int W_{\theta_i}(x) \widehat\rho(\omega,\mathrm{d}x)
\widehat
Q_{\theta_i}(\mathrm{d}\omega)=0$ for $i = 1,2$. On the other hand, any Markov
kernel $\rho\dvtx[0,\infty]^4\times\BBB(D)\to[0,1]$ is of the form
$\rho(\omega,B)=p(\omega)\epsilon_{(\beta_1,\alpha_1)}(B)+
(1-p(\omega))\epsilon_{(\beta_2,\alpha_2)}(B)$, where
$p\dvtx[0,\infty]^4\to[0,1]$ is Borel, $\omega\in[0,\infty]^4$ and
$B\in\BBB(D)$. It is easy to see that for such a Markov kernel $\rho$,
there exists a Markov kernel $ \bar\rho\dvtx[0,\infty]^4\times
\BBB(D)\to[0,1]$ such that both
\[
\int W_{\theta_1}(x) \rho(\omega,\mathrm{d}x) Q_{\theta_1}(\mathrm{d}\omega)\ge
\mathrm{e}^{-\gamma} \bigl(1-p(h_0,\infty,\infty,\infty)\bigr)
\]
and
\[
\hspace*{-22pt}\int W_{\theta_2}(x) \rho(\omega,\mathrm{d}x) Q_{\theta_2}(\mathrm{d}\omega)\ge
\mathrm{e}^{-\gamma}p(h_0,\infty,\infty,\infty) .
\]
In view of
\eqref{defdeffic}, we thus have
$\delta(\FFF(\Theta),\widehat\FFF(\Theta))\ge
\delta(\FFF(\Theta_0),\FFF(\Theta_0))\ge \mathrm{e}^{-\gamma}/2$,
which completes the proof of (ii).

(iii)   This follows by the same arguments as in (ii).

\setcounter{equation}{0} \setcounter{subsection}{1}

\begin{appendix}\label{appendix}
\section*{Appendix}

Here, we collect necessary facts regarding Le Cam's distance in
deficiency. The reader is referred to Le Cam \cite{LC86}, Le Cam and
Young \cite{LY90} and Strasser's monograph \cite{St85} for unexplained
notation encountered in this section. Let $\Theta$ be a non-empty set,
$(E,\AAA)$ be a measurable space and $(P_{\theta})_{\theta\in\Theta}$
be a family of probability measures on $\AAA$. The triplet
$\EEE=(E,\AAA,(P_{\theta})_{\theta\in\Theta})$ is then called a~(statistical) experiment. Consider two experiments
$\EEE_i=(E_i,\AAA_i,(P_{i,\theta})_{\theta\in\Theta})$, $i=1,2$,
indexed by $\Theta$. A decision problem is a triple $(\Theta,D,W)$,
where $D$ is a topological space and $W=(W_\theta)_{\theta\in\Theta}$
is a loss function $W_\theta\dvtx D\to\RR$, $\theta\in\Theta$. Let
$\|W_\theta\|_{\infty}=\sup_{d\in D}|W_\theta(d)|$. Also, let
$\epsilon\ge0$. Then, $\EEE_1$ is called $\epsilon$-\textit{deficient
with respect to} $\EEE_2$, notated as $\EEE_1\supseteq_\epsilon
\EEE_2$, if and only if for all decision problems $(\Theta,D,W)$ with
$W$ continuous and bounded, and all $\beta_2\in\BBB(\EEE_2,D)$, there
exists some $\beta_1\in\BBB(\EEE_1,D)$ such that\looseness=-1
\[
\beta_1(W_\theta,P_{1,\theta})\le\beta_2(W_\theta,P_{2,\theta
})+\epsilon\|W_\theta\|_{\infty} ,\qquad  \theta\in\Theta ,
\]\looseness=0
where $\BBB(\EEE_i,D)$ ($i=1,2$) is the space of generalized
decision functions (see \cite{St85}, Definition~42.2). The deficiency
of $\EEE_1$ with respect to $\EEE_2$ is the number
%
\begin{equation}\label{defdeffic}
\delta(\EEE_1,\EEE_2) = \inf\{\epsilon>0 \dvtx \EEE_1\supseteq
_\epsilon
\EEE_2\} .
\end{equation}
The relation $\EEE_1\supseteq_\epsilon
\EEE_2$ is interpreted in the following sense: we have
$\EEE_1\supseteq_\epsilon\EEE_2$ if $\EEE_1$ is more informative than
$\EEE_2$ uniformly over all decision problems with continuous and
bounded loss functions up to some error $\epsilon$. Two experiments
$\EEE_1$ and $\EEE_2$ are called \textit{equivalent in deficiency} if and
only if $\EEE_1\supseteq_0 \EEE_2$ and $\EEE_2\supseteq_0 \EEE_1$.

Recall that (see \cite{St85}, Lemma 55.4 and Remark 55.6(2))
%
\begin{equation}\label{defMarkov}
\delta(\EEE_1,\EEE_2) = \inf_K \sup_{\theta\in\Theta}  \|
P_{2,\theta} - K P_{1,\theta} \|
\end{equation}
with infimum now taken over all Markov kernels $K\dvtx E_1\times\EEE_2\to
[0,1]$.

Le Cam's distance between $\EEE_1$ and $\EEE_2$ is a pseudo-metric on
the space of all experiments indexed by $\Theta$ (see
\cite{St85}, Corollary 59.6), defined by setting
%
\begin{equation}\label{deflecam}
\Delta(\EEE_1,\EEE_2) = \max\{
\delta(\EEE_1,\EEE_2),\delta(\EEE_2,\EEE_1) \}.
\end{equation}
If $(E_1,\AAA_1)=(E_2,\AAA_2)$, then we have (see
\cite{St85}, Corollary 59.6)
%
\begin{equation}\label{deflecamsmallerTV}
\Delta(\EEE_1,\EEE_2) \le
\sup_{\theta\in\Theta}\|P_{1,\theta}-P_{2,\theta}\| .
\end{equation}
Clearly, if $\EEE_1$ and $\EEE_2$ are two experiments indexed by the
same $\Theta$, then $\EEE_1$ is equivalent to $\EEE_2$ in deficiency if
and only if $\Delta(\EEE_1,\EEE_2)=0$. Let $\EEE$, $\EEE_n$, $\FFF_n$,
$n\in\NN$, be experiments, all indexed by $\Theta$. We then say that
$\EEE_n$ converges (strongly) in deficiency, or that $\EEE_n$ and~$\FFF_n$ are (strongly) asymptotically equivalent in deficiency, if and
only if $\Delta(\EEE_n,\EEE) \to0$ and $\Delta(\EEE_n,\FFF_n)
\to0$, respectively, as $n\to\infty$.

For $\varnothing\neq\Theta_0\subseteq\Theta$, we employ the notation
$\EEE(\Theta_0)=(E,\AAA,(P_{\theta})_{\theta_\in\Theta_0})$ for
corresponding subexperiments of
$\EEE=(E,\AAA,(P_\theta)_{\theta\in\Theta})$. We refer to weak
convergence and weak asymptotic equivalence in deficiency if and only
if, for all non-empty and finite $\Theta_0\subseteq\Theta$, the
corresponding subexperiments converge strongly and are strongly
asymptotically equivalent in deficiency, respectively.
\end{appendix}

\section*{Acknowledgements}
This research was supported in part by ARC Grant DP0988483 and DFG
Grant 447 AUS-111/3/07.

\printhistory


\begin{thebibliography}{28}

\bibitem{AMZh}
A{\"{i}}t-Sahalia, Y., Mykland, P.A. and Zhang, L. (2011).
Ultra high frequency volatility estimation with dependent
microstructure noise.
\textit{J. Econometrics} \textbf{160} 160--175.

\bibitem{Bi68}
Billingsley, P. (1968).\vadjust{\goodbreak} \textit{Convergence of Probability
Measures.} New York: Wiley.
\MR{0233396}

\bibitem{BS73}
Black, F. and Scholes, M. (1973). The pricing of options and
corporate liabilities. \textit{J. Political Economy} \textbf{81} 637--659.

\bibitem{Bo86}
Bollerslev, T. (1986). Generalized autoregressive
conditional heteroscedasticity. \textit{J. Econometrics} \textbf{31} 307--327.
\MR{0853051}

\bibitem{BL96}
Brown, L.D. and Low, M. (1996). Asymptotic equivalence of
nonparametric regression and white noise. \textit{Ann. Statist.}
\textbf{24} 2384--2398.
\MR{1425958}

\bibitem{BWZ03}
Brown, L.D., Wang, Y. and Zhao, L.H. (2003). On the
statistical equivalence at suitable frequencies of GARCH and
stochastic volatility models with the corresponding diffusion model.
\textit{Statist. Sinica} \textbf{13} 993--1013.
\MR{2026059}

\bibitem{C07}
Carter, A.V. (2007). Asymptotic approximation of nonparametric
regression experiments with unknown variances.
\textit{Ann. Statist.} \textbf{35} 1644--1673.
\MR{2351100}

\bibitem{DR06}
Dalalyan, A.S. and Reiss, M. (2006). Asymptotic
statistical equivalence for scalar ergodic diffusions.
\textit{Probab. Theory Related Fields} \textbf{134} 248--282.
\MR{2222384}

\bibitem{EK97}
Embrechts, P., Kl{\"{u}}ppelberg, C. and Mikosch, T. (1997).
\textit{Modelling Extremal Events for Insurance and Finance.}
Berlin: Springer.
\MR{1458613}

\bibitem{Engle82}
Engle, R.F. (1982). Autoregressive conditional
heteroscedasticity with estimates of the United Kingdom inflation.
\textit{Econometrica} \textbf{50} 987--1007.
\MR{0666121}

\bibitem{FHR00}
Falk, M., H\"{u}sler, J. and Reiss, R.D. (1994).
\textit{Laws of Small Numbers: Extremes and Rare Events.}
Basel: Birkh\"{a}user.
\MR{1296464}

\bibitem{FKL06}
Fasen, V., Kl\"{u}ppelberg, C. and Lindner, A. (2006). Extremal
behavior of stochastic volatility models. In \textit{Stochastic
Finance} (A. Shiryaev,
M.D.R. Grossinho, P. Oliviera, M. Esquivel, eds.) 107--155. New York:
Springer.
\MR{2230762}

\bibitem{GN06}
Grama, I.O. and Neumann, M.H. (2006). Asymptotic equivalence
of nonparametric autoregression and nonparametric regression.
\textit{Ann. Statist.} \textbf{34} 1701--1732.
\MR{2283714}

\bibitem{HubPos}
Hubalek, F. and Posedel, P. (2010).
Joint analysis and estimation of stock prices and trading volume
in Barndorff-Nielsen and Shephard stochastic volatility models.
\textit{Quant. Finance}. To appear.

\bibitem{JaKlMu}
Jacod, J., Kl\"{u}ppelberg, C. and M\"{u}ller, G. (2010).
Testing for non-correlation or for a~functional relationship
between price and volatility jumps.
Preprint, Univ. Paris VI and Technische Univ. M\"{u}nchen.

\bibitem{KV07}
Kallsen, J. and Vesenmayer, B. (2009). COGARCH as a
continuous time limit of GARCH(1,1). \textit{Stochastic
Process. Appl.} \textbf{119} 74--98.
\MR{2485020}

\bibitem{KLM04}
Kl{\"{u}}ppelberg, C., Lindner, A. and Maller,
R. (2004). A continuous-time {GARCH} process driven by a {L}\'{e}vy
process: Stationarity and second-order behaviour. \textit{J. Appl.
Probab.} \textbf{41} 601--622.
\MR{2074811}

\bibitem{LC86}
Le Cam, L. (1986). \textit{Asymptotic Methods in Statistical
Decision Theory.} New York: Springer.
\MR{0856411}

\bibitem{LY90}
Le Cam, L. and Yang, G.L. (1990). \textit{Asymptotics in
Statistics. Some Basic Concepts.} New York: Springer.
\MR{1066869}

\bibitem{MMS08}
Maller, R., M\"{u}ller, G. and Szimayer, A. (2008). GARCH
modelling in continuous time for irregularly spaced time series data.
\textit{Bernoulli} \textbf{14} 519--542.
\MR{2544100}

\bibitem{Me73}
Merton, R. (1973). The theory of rational option pricing.
\textit{Bell J. Econom. Management Sci.} \textbf{4} 141--183.
\MR{0496534}

\bibitem{MN98}
Milstein, G. and Nussbaum, M. (1998). Diffusion
approximation for nonparametric autoregression. \textit{Probab. Theory
Related Fields} \textbf{112} 535--543.
\MR{1664703}

\bibitem{Ne90}
Nelson, D.B. (1990). ARCH models as diffusion
approximations. \textit{J. Econometrics} \textbf{45} 7--38.
\MR{1067229}

\bibitem{Nu96}
Nussbaum, M. (1996). Asymtotic equivalence of density
estimation and Gaussian white noise. \textit{Ann. Statist.} \textbf{24}
2399--2430.
\MR{1425959}

\bibitem{Re87}
Resnick, S.I. (1987). \textit{Extreme Values, Regular Variation, and
Point Processes.} New York: Springer.
\MR{0900810}

\bibitem{Re93}
Reiss, R.-D. (1993). \textit{A Course on Point Processes.}
New York: Springer.
\MR{1199815}

\bibitem{St85}
Strasser, H. (1985). \textit{Mathematical Theory of Statistics.}
Berlin: de Gruyter.
\MR{0812467}

\bibitem{Wa02}
Wang, Y. (2002). Asymptotic nonequivalence of GARCH models
and diffusions. \textit{Ann. Statist.} \textbf{30} 754--783.
\MR{1922541}

\end{thebibliography}
\end{document}